\documentclass{amsart}

\usepackage{amsthm}
\usepackage{amsmath}
\usepackage{amssymb}
\usepackage{amsfonts}
\usepackage{latexsym}
\usepackage[all]{xy} \SelectTips{eu}{}
\usepackage{hyperref}
\usepackage{eucal}
\usepackage{xcolor}
\usepackage{marginnote}
\usepackage[nameinlink,capitalise,noabbrev]{cleveref}

\hypersetup{
     colorlinks=true,
     linkcolor=blue!60!black,
     filecolor=blue!75!black,
     citecolor = green!50!black,
     urlcolor=magenta,
     }

\newtheorem{intthm}{Theorem}[]

\newtheorem*{intque*}{Question}
\newtheorem*{intexa*}{Example}
\newcommand{\numberseries}{\bfseries}   
\newlength{\thmtopspace}                
\newlength{\thmbotspace}                
\newlength{\thmheadspace}               
\newlength{\thmindent}                  
\newtheoremstyle{bfupright head,slanted body}
{\thmtopspace}{\thmbotspace}
{\slshape}{\thmindent}{\bfseries}{.}{\thmheadspace}
{{\numberseries \thmnumber{#2\;}}\thmnote{#3}}

\newtheoremstyle{bfupright head,upright body}
{\thmtopspace}{\thmbotspace}
{\upshape}{\thmindent}{\bfseries}{.}{\thmheadspace}
{{\numberseries \thmnumber{#2\;}}\thmnote{#3}}

\newtheoremstyle{fixed bf head,slanted body}
{\thmtopspace}{\thmbotspace}{\slshape}
{\thmindent}{\bfseries}{.}{\thmheadspace}
{{\numberseries \thmnumber{#2\;}}\thmname{#1}\thmnote{ (#3)}}

\newtheoremstyle{fixed bf head,upright body}
{\thmtopspace}{\thmbotspace}{\upshape}
{\thmindent}{\bfseries}{.}{\thmheadspace}
{{\numberseries \thmnumber{#2\;}}\thmname{#1}\thmnote{ (#3)}}

\newtheoremstyle{numbered paragraph}
{\thmtopspace}{\thmbotspace}{\upshape}
{\thmindent}{\upshape}{}{\thmheadspace}
{{\numberseries \thmnumber{#2.}}}

\theoremstyle{bfupright head,slanted body}
\newtheorem{res}{}[section]             \newtheorem*{res*}{}

\theoremstyle{bfupright head,upright body}
\newtheorem{bfhpg}[res]{}               \newtheorem*{bfhpg*}{}

\theoremstyle{fixed bf head,slanted body}
\newtheorem{thm}[res]{Theorem}          \newtheorem*{thm*}{Theorem}
\newtheorem{prp}[res]{Proposition}      \newtheorem*{prp*}{Proposition}
\newtheorem{cor}[res]{Corollary}        \newtheorem*{cor*}{Corollary}
\newtheorem{lem}[res]{Lemma}            \newtheorem*{lem*}{Lemma}
         \newtheorem*{que*}{Question}
\theoremstyle{fixed bf head,upright body}
\newtheorem{dfn}[res]{Definition}       \newtheorem*{dfn*}{Definition}
\newtheorem{rmk}[res]{Remark}           \newtheorem*{rmk*}{Remark}
\newtheorem{exa}[res]{Example}           \newtheorem*{exa*}{Example}
\newtheorem{setup}[res]{Setup}           \newtheorem*{setup*}{Setup}
           \newtheorem*{nota*}{Notation}

\theoremstyle{numbered paragraph}
\newtheorem{ipg}[res]{}

\newlength{\thmlistleft}        
\newlength{\thmlistright}       
\newlength{\thmlistpartopsep}   
\newlength{\thmlisttopsep}      
\newlength{\thmlistparsep}      
\newlength{\thmlistitemsep}     

\newcounter{eqc}
\newenvironment{eqc}{\begin{list}{\upshape (\textit{\roman{eqc}})}%
		{\usecounter{eqc}%
			\setlength{\leftmargin}{\thmlistleft}%
			\setlength{\labelwidth}{\thmlistleft}%
			\setlength{\rightmargin}{\thmlistright}%
			\setlength{\partopsep}{\thmlistpartopsep}%
			\setlength{\topsep}{\thmlisttopsep}%
			\setlength{\parsep}{\thmlistparsep}%
			\setlength{\itemsep}{\thmlistitemsep}}}%
	{\end{list}}%

\newcounter{prt}
\newenvironment{prt}{\begin{list}{\upshape (\alph{prt})}%
		{\usecounter{prt}%
			\setlength{\leftmargin}{\thmlistleft}%
			\setlength{\labelwidth}{\thmlistleft}%
			\setlength{\rightmargin}{\thmlistright}%
			\setlength{\partopsep}{\thmlistpartopsep}%
			\setlength{\topsep}{\thmlisttopsep}%
			\setlength{\parsep}{\thmlistparsep}%
			\setlength{\itemsep}{\thmlistitemsep}}}%
	{\end{list}}%

\newcounter{rqm}
\newenvironment{rqm}{\begin{list}{\upshape (\arabic{rqm})}%
		{\usecounter{rqm}%
			\setlength{\leftmargin}{\thmlistleft}%
			\setlength{\labelwidth}{\thmlistleft}%
			\setlength{\rightmargin}{\thmlistright}%
			\setlength{\partopsep}{\thmlistpartopsep}%
			\setlength{\topsep}{\thmlisttopsep}%
			\setlength{\parsep}{\thmlistparsep}%
			\setlength{\itemsep}{\thmlistitemsep}}}%
	{\end{list}}%

	{\end{list}}%

\newenvironment{prf*}[1][Proof]{%
	\begin{proof}[\bf #1]
		\setcounter{equation}{0}
		}
	{\end{proof}
}

\newcommand{\pgref}[1]{\ref{#1}}

\renewcommand{\eqref}[1]{(\pgref{eq:#1})}

\newcommand{\eqclbl}[1]{{\upshape(\textit{#1})}}
\newcommand{\proofofimp}[3][:]{\mbox{\eqclbl{#2}$\!\implies\!$\eqclbl{#3}#1}}
\numberwithin{equation}{res}
\def\urltilda{\kern -.15em\lower .7ex\hbox{\~{}}\kern .04em}

\newcommand{\setof}[3][\mspace{1mu}]{\{#1#2 \mid #3#1\}}
\newcommand{\GP}[1]{\mathsf{GP}(#1)}
\newcommand{\TacQ}{{}_{\Q,A}\mathsf{Tac}(\sf Prj)}
\newcommand{\AcQ}{{}_{\Q,A}\mathsf{Acy}(\sf Prj)}
\newcommand{\Mod}[1]{{_{#1}\mathsf{Mod}}}
\newcommand{\Prj}[1]{\mathsf{Prj}(#1)}
\newcommand{\Inj}[1]{\mathsf{Inj}(#1)}
\newcommand{\ch}[1]{\mathsf{Ch}(#1)}

\newcommand{\xra}[2][]{\xrightarrow[#1]{\:#2\:}}

\newcommand{\Ker}[1]{\nobreak{\operatorname{Ker}#1}}
\newcommand{\Hom}{\operatorname{Hom}}
\newcommand{\Ext}[4][]{\operatorname{Ext}_{#1}^{#2}(#3,#4)}
\newcommand{\Tor}[4][]{\operatorname{Tor}_{#2}^{#1}(#3,#4)}
\newcommand{\Coker}[1]{\nobreak{\operatorname{Coker}#1}}
\newcommand{\C}{\mathcal{C}}
\newcommand{\Q}{\mathsf{Q}}
\newcommand{\Qop}{\Q^{\sf op}}
\newcommand{\pd}{\mbox{\rm pd}}

\newcommand{\id}{\mathrm{id}}
\newcommand{\kk}{\Bbbk}
\newcommand{\is}{\cong}

\newcommand{\lMod}[1]{{}_{#1}\mspace{-1mu}{\sf{Mod}}}
\newcommand{\lE}[1]{{}_{#1}\mspace{-1mu}{\sf{Acy}}}

\newcommand{\lPrj}[1]{{}_{#1}\mspace{-1mu}{\sf{Prj}}}
\newcommand{\rMod}[1]{{\sf{Mod}}_{#1}}
\newcommand{\Aop}{A^{\mathsf{op}}}
\newcommand{\nablaq}{\nabla_q}
\newcommand{\DeltaLq}{{}_q\Delta}
\newcommand{\DeltaRq}{\Delta_q}
\def\soft#1{\leavevmode\setbox0=\hbox{h}\dimen7=\ht0\advance
	\dimen7 by-1ex\relax\if t#1\relax\rlap{\raise.6\dimen7
		\hbox{\kern.3ex\char'47}}#1\relax\else\if T#1\relax
	\rlap{\raise.5\dimen7\hbox{\kern1.3ex\char'47}}#1\relax
	\else\if d#1\relax\rlap{\raise.5\dimen7\hbox{\kern.9ex
			\char'47}}#1\relax\else\if D#1\relax\rlap{\raise.5\dimen7
		\hbox{\kern1.4ex\char'47}}#1\relax\else\if l#1\relax
	\rlap{\raise.5\dimen7\hbox{\kern.4ex\char'47}}#1\relax
	\else\if L#1\relax\rlap{\raise.5\dimen7\hbox{\kern.7ex
			\char'47}}#1\relax\else\message{accent \string\soft
		\space #1 not defined!}#1\relax\fi\fi\fi\fi\fi\fi}

\makeatletter
\@namedef{subjclassname@2020}{%
  \textup{2020} Mathematics Subject Classification}
\makeatother

\begin{document}

\title[Total acyclicity of Q-shaped diagrams and N-stable categories]
{Total acyclicity of Q-shaped diagrams and N-stable categories}

\author[L. Liang]{Li Liang}

\address{Li Liang: Department of Mathematics, Gansu Provincial Basic Discipline Research Center for Complex Systems Analysis and Control, Lanzhou Jiaotong University, Lanzhou 730070, China}

\email{lliangnju@gmail.com}

\urladdr{https://sites.google.com/site/lliangnju}

\thanks{L. Liang was partly supported by NSF of China (Grant No. 12271230) and NSF of Gansu Province (Grant No. 26JRRA059).}

\keywords{$\Q$-shaped diagram, Gorenstein projective module, defect category, $N$-stable category.}

\subjclass[2020]{18G25, 18A25, 18G35.}

\begin{abstract}
We study total acyclicity of $\Q$-shaped diagrams of projective modules. We prove that an acyclic $\Q$-shaped diagram $P$ of projective $A$-modules is totally acyclic if and only if the associated module $K_q(P)$, or equivalently $C_q(P)$, is Gorenstein projective for every $q\in \Q$. We also establish Frobenius structures on the categories of acyclic and totally acyclic $\Q$-shaped diagrams of projective $A$-modules and introduce the associated defect category. Using adjunctions between $\Q$-shaped diagrams and chain complexes, we compare this category with the classical projective defect category and show that their vanishing is equivalent under a mild condition on the Serre functor of $\Q$. Finally, we apply our results to $N$-complexes and extend the $N$-stable equivalence of Brightbill and Miemietz beyond the Gorenstein setting, obtaining a triangle equivalence between the homotopy category of totally $N$-acyclic complexes of projective $A$-modules and the $N$-stable category of Gorenstein projective $A$-modules.
\end{abstract}

\maketitle
\vspace*{.3cm}
\vspace*{-.9cm}
\enlargethispage{.4cm}

\section*{Introduction}
\noindent
{\bf Throughout the paper}, $\kk$ is a commutative noetherian hereditary ring (e.g. $\kk$ is a field or $\kk=\mathbb{Z}$), $A$ is an arbitrary $\kk$-algebra, and $\Q$ is a small $\kk$-linear category satisfying the hypotheses of \cite[Setup 2.9]{HJ24}; see Setup \ref{setup:Q}. By an $A$-module we mean a left $A$-module; right $A$-modules are modules over $\Aop$. We write $\lMod{\Q,A}$ for the category of $\kk$-linear functors from $\Q$ to the category $\lMod{A}$ of $A$-modules.

The distinction between acyclicity and total acyclicity is a basic theme in Gorenstein homological algebra. Recall that an acyclic complex $X$ of projective $A$-modules is \emph{totally acyclic} if $\Hom_A(X,L)$ is acyclic for every projective $A$-module $L$. The cycles of such complexes are precisely the Gorenstein projective $A$-modules, so total acyclicity provides a homological link between complexes and Gorenstein projective objects; see e.g. \cite{EnochsGP,HHl04a}. We mention that the distinction can also be expressed at the level of triangulated categories. Writing $K_{\mathrm{ac}}(\Prj A)$ and $K_{\mathrm{tac}}(\Prj A)$ for the homotopy categories of acyclic and totally acyclic complexes of projective $A$-modules, respectively, one gets the Verdier quotient
$$\mathsf D_{\mathrm{def}}(A)=K_{\mathrm{ac}}(\Prj A)\big/K_{\mathrm{tac}}(\Prj A),$$
which records the difference between the two notions. This Verdier quotient, known as the projective defect category, was studied by Iyengar and Krause \cite{IyengarKrause06}, who investigated its structure and its relationship with Gorenstein properties of commutative noetherian rings admitting dualizing complexes.

The theory of $\Q$-shaped derived categories extends the homological algebra of complexes to diagrams indexed by suitable small categories. Indeed, a chain complex can be viewed as a representation of the infinite linear quiver subject to the relations that two consecutive arrows compose to zero. Motivated by the role of a Serre functor in this description, Holm and J{\o}rgensen developed the $\Q$-shaped derived category and its homological framework in
\cite{HJ22,HJ24,HJ24b}. The resulting setting includes, among other examples, $N$-complexes, for which the composite of any $N$ consecutive differentials vanishes. In this generality, the notion of an acyclic (or exact) $\Q$-shaped diagram is formulated in terms of the homological
structure of $\Q$, rather than ordinary homology groups.

More recently, Holm and J{\o}rgensen \cite{HJ26} introduced and studied totally acyclic $\Q$-shaped diagrams of projective modules. Their work establishes an adjoint triple relating $\Q$-shaped diagrams to chain complexes and, for commutative noetherian coefficient rings, characterizes the Gorenstein property by the coincidence of acyclicity and total acyclicity under a mild condition on the Serre functor; see \cite[Theorem B]{HJ26}. Their results suggest two related questions: can the difference between acyclicity and total acyclicity be compared with its classical counterpart at the level of triangulated categories, and can total acyclicity of an individual $\Q$-shaped diagram be detected by Gorenstein projective $A$-modules? We address both questions without requiring the coefficient ring to be Gorenstein.

We begin with the triangulated aspect of this problem. Let $\AcQ$ and $\TacQ$ denote the categories of acyclic and totally acyclic $\Q$-shaped \mbox{diagrams} of projective $A$-modules, respectively. We show that they are Frobenius exact categories and the stable category $\mathsf{stab}(\TacQ)$ is a thick triangulated subcategory of $\mathsf{stab}(\AcQ)$; see Propositions \ref{thm:EPrj-Frobenius} and \ref{thm:Tac-Frobenius}. This gives the \emph{$\Q$-shaped projective defect category}
$$\mathsf D_{\mathrm{def}}^{\Q}(A)=\mathsf{stab}(\AcQ)\big/\mathsf{stab}(\TacQ),$$
which is the natural counterpart of $\mathsf D_{\mathrm{def}}(A)$. We prove that, for each $q\in\Q$, the adjoint pair $({}_q\Delta, \nabla_q)$ from \cite{HJ26} induces an adjoint pair of triangle functors
$${}_q\overline{\Delta}:\mathsf D_{\mathrm{def}}(A) \rightleftarrows
\mathsf D_{\mathrm{def}}^{\Q}(A):\overline{\nabla}_q;$$
see Proposition \ref{thm:defect-adjunction}. The family $\{\overline{\nabla}_q\}_{q\in\Q}$ is jointly conservative, while ${}_q\overline{\Delta}$ is conservative whenever $\mathbb S q\neq q$, where $\mathbb S$ is the Serre functor of $\Q$; see Propositions \ref{cor:defect-detection} and \ref{thm:defect-left-conservative}. These adjunctions then lead to our first main result, which compares the vanishing of the two defect categories; see Theorem \ref{thm:classical-Q-defect-equivalence}.

\begin{intthm}\label{thmA}
Consider the following conditions:
\begin{eqc}
\item $K_{\mathrm{ac}}(\Prj A)=K_{\mathrm{tac}}(\Prj A)$.
\item $\AcQ=\TacQ$.
\end{eqc}
Then \eqclbl{i}$\Longrightarrow$\eqclbl{ii} holds. If, in addition, there exists $q\in\Q$ such that $\mathbb S q\neq q$, then the two conditions are equivalent.
\end{intthm}

We mention that the hypothesis in the implication \eqclbl{ii}$\Longrightarrow$\eqclbl{i} cannot be omitted; see Example \ref{exa:converse-fails}. When $A$ is commutative noetherian, the above theorem, together with the classical criterion, recovers the ring-theoretic characterization of Holm and J{\o}rgensen \cite[Theorem B]{HJ26}. More generally, the result identifies when the
$\Q$-shaped and classical projective defect categories vanish simultaneously, even over noncommutative coefficient rings $A$.

The comparison of defect categories is global; it concerns when every acyclic $\Q$-shaped diagram of projective $A$-modules is totally acyclic. To understand the objects entering this comparison, it is natural to seek an analogue of the classical criterion in terms of cycles. Unlike an ordinary complex, a general $\Q$-shaped diagram has no single family of cycle modules. Instead, for each $q\in\Q$, the stalk functor $S_q$ is part of an adjoint triple $(C_q, S_q, K_q)$, with $C_q$ and $K_q$ taking values in $A$-modules; see \ref{stalk}. These modules provide the appropriate tests for total acyclicity of $\Q$-shaped diagrams. Our second main result makes this
precise; see Theorem \ref{thm:totally-acyclic-characterization}.

\begin{intthm}\label{thmB}
Let $P$ be an acyclic object in $\lMod{\Q,A}$ with $P(q)$ projective for every $q\in\Q$. Then the following statements are equivalent.
\begin{eqc}
\item $P$ is a totally acyclic $\Q$-shaped diagram of projective
$A$-modules.
\item $K_q(P)$ is Gorenstein projective for every $q\in\Q$.
\item $C_q(P)$ is Gorenstein projective for every $q\in\Q$.
\end{eqc}
\end{intthm}

Thus, the total acyclicity of an acyclic $\Q$-shaped diagram $P$ of projective modules is characterized entirely in terms of the Gorenstein projectivity of the modules $K_q(P)$ or $C_q(P)$. No commutativity or Gorenstein assumption on $A$ is needed.

\begin{equation*}
  \ast \ \ \ast \ \ \ast
\end{equation*}
\noindent
We apply our results to $N$-complexes, an important instance of the $\Q$-shaped setting, whose homological theory was developed, in particular, by Iyama, Kato and Miyachi \cite{IKM17}. In this context, Brightbill and Miemietz \cite[Theorem 4.12]{BM24} established an $N$-analogue of Buchweitz's equivalence for Gorenstein abelian categories, identifying the homotopy category of acyclic $N$-complexes of projectives with an $N$-stable category.

Our preceding results allow us to extend this equivalence beyond the Gorenstein setting. The key ingredient is our characterization of total acyclicity specialized to $N$-complexes, that is, an $N$-acyclic $N$-complex $P$ of projective $A$-modules is totally $N$-acyclic if and only if all its generalized cycles $Z_n^r(P)$ are Gorenstein projective for $n\in\mathbb Z$ and $1\leq r<N$; see Theorem \ref{cor:N-cycles-GP}. Consequently, the generalized cycles of a totally $N$-acyclic $N$-complex of projectives give rise to a chain of admissible monomorphisms in $\GP{A}$, the category of Gorenstein projective $A$-modules. Conversely, adapting the construction of Brightbill and Miemietz \cite[Proposition 4.5]{BM24}, such a chain can be
extended to an $N$-acyclic $N$-complex of projectives whose generalized cycles are Gorenstein projective; the resulting complex is totally $N$-acyclic by Theorem \ref{cor:N-cycles-GP}. Together with the Frobenius structures established earlier, this yields a triangle equivalence with the $N$-stable category $\mathsf{stab}^N(\GP{A})$ of chains of $N-2$ admissible monomorphisms in $\GP{A}$. Our third main result is as follows; see Theorem \ref{thm:N-stable-total}.

\begin{intthm}\label{thmC}
For every $N\geq2$, there is a triangle equivalence
$$K^N_{\mathrm{tac}}(\Prj A)\simeq
\mathsf{stab}^N(\GP{A}).$$
In particular, for $N=2$, this recovers the classical equivalence
$$K_{\mathrm{tac}}(\Prj A)\simeq\mathsf{stab}(\GP{A}).$$
\end{intthm}

Here, $K^N_{\mathrm{ac}}(\Prj A)$ and $K^N_{\mathrm{tac}}(\Prj A)$ are the homotopy categories of $N$-acyclic and totally $N$-acyclic $N$-complexes of projective $A$-modules, respectively. If $A$ is Iwanaga--Gorenstein, then every $N$-acyclic $N$-complex of projective $A$-modules is totally $N$-acyclic, and so Theorem \ref{thmC} recovers \cite[Theorem 4.12]{BM24}.

\section{Preliminaries}\label{property}
\noindent
In this section, we record some basic notions and facts which are used in the paper.

\begin{setup}\label{setup:Q}
Throughout the paper, we assume that $\Q$ is a small $\kk$-linear
category satisfying the hypotheses of \cite[Setup 2.9]{HJ24}:
\begin{rqm}
\item \emph{Hom-finiteness}: Each $\kk$-module $\Q(p,q)$ is finitely generated and projective.
\item \emph{Local Boundedness}: For every object $q$ in $\Q$, there are only finitely many objects in $\Q$ mapping nontrivially into or out of $q$.
\item \emph{Existence of a Serre Functor}: There is a $\kk$-linear auto-equivalence $\mathbb{S}: \Q \to \Q$ such that $\Q(p,q)\is\Hom_{\kk}(\Q(q,\mathbb{S}(p)),\kk)$ for all $p,q\in\Q$.
\item \emph{Strong Retraction Property}: For each object $q\in\Q$, the unit map $\kk\to \Q(q,q)$, given by $x\mapsto x\cdot\id_q$, admits a $\kk$-module retraction and there exists a collection $\{\mathfrak{r}_q\}_{q\in\Q}$ of complements, that is, $\kk$-modules $\mathfrak{r}_q$ such that
      $$\Q(q,q)=(\kk\cdot\id_q)\oplus\mathfrak{r}_q,$$
      compatible with composition in $\Q$ in the sense that
    \begin{eqc}
    \item[($\dagger$)] $\mathfrak{r}_q\circ\mathfrak{r}_q\subseteq\mathfrak{r}_q$ for all $q\in\Q$;
    \item[($\ddagger$)] $\Q(q,p)\circ\Q(p,q)\subseteq\mathfrak{r}_q$ for all $p\neq q$ in $\Q$.
    \end{eqc}
\item The pseudo-radical $\mathfrak{r}$ \footnote{The Strong Retraction Property makes it possible to define the pseudo-radical $\mathfrak{r}$; see \cite[Remark 7.4 and Lemma 7.7]{HJ22}.} is nilpotent, that is, $\mathfrak{r}^N=0$ for some $N\in\mathbb{N}$.
\end{rqm}

{\bf Throughout this paper}, we denote by $\lMod{\Q,A}$ the category of $\kk$-linear functors from $\Q$ to the (left) $A$-module category $\lMod{A}$; such functors are also called the $\Q$-shaped diagrams of $A$-modules. The hom set ($\kk$-module) in $\lMod{\Q,A}$ is written $\Hom_{\Q,A}$ and its right derived functors are denoted by $\mathrm{Ext}_{\Q,A}^i$. We also denote by $\rMod{\Q,A}$ the category of $\kk$-linear functors from $\Q^{\sf op}$ to the $\Aop$-module category $\rMod{A}$. In the case where $A=\kk$, we write $\lMod{\Q}$ (resp., $\rMod{\Q}$) instead of $\lMod{\Q,A}$ (resp., $\rMod{\Q,A}$).
\end{setup}

\begin{bfhpg}[Evaluation functor and adjunctions]\label{eva}
Following \cite[Corollary 3.9]{HJ22}, for each $q \in \Q$, there is an adjoint triple $(F_{q},E_{q},G_{q})$ defined as follows:
  \begin{equation*}
  \xymatrix@C=4pc{
      \lMod{\Q,A}
    \ar[r]^-{E_{q}}
    & \lMod{A}
    \ar@/_1.8pc/[l]_-{F_{q}}
    \ar@/^1.8pc/[l]^-{G_{q}}
  }
  \qquad \text{given by} \qquad
  {\setlength\arraycolsep{1.5pt}
   \renewcommand{\arraystretch}{1.2}
  \begin{array}{rcl}
  F_q(M) &=& \Q(q,-) \otimes_\kk M \\
  E_q(X) &=& X\mspace{1.5mu}(q) \\
  G_{q}(M) &=& \Hom_{\kk}(\Q(-,q),M)\;.
  \end{array}
  }
  \end{equation*}
Note that all the functors $F_q$, $E_q$ and $G_q$ are exact as $\Q$ satisfies Hom-finiteness.
\end{bfhpg}

\begin{bfhpg}[Stalk functor and adjunctions]\label{stalk}
Following \cite[Proposition 7.15]{HJ22}, for each $q\in\Q$, there is an adjoint triple $(C_{q},S_{q},K_{q})$ defined as follows:
  \begin{equation*}
  \xymatrix@C=4pc{
    \Mod{A}
    \ar[r]^-{S_{q}}
    &
    \Mod{\C,A}
    \ar@/_1.8pc/[l]_-{C_{q}}
    \ar@/^1.8pc/[l]^-{K_{q}}
  }
  \qquad \text{given by} \qquad
  {\setlength\arraycolsep{1.5pt}
   \renewcommand{\arraystretch}{1.2}
  \begin{array}{rcl}
  C_q(X) &=& S\{q\} \otimes_\Q X \\
  S_q(M) &=& S\langle q\rangle\otimes_\kk M \\
  K_{q}(X) &=& \Hom_\Q(S\langle q\rangle,X)\;.
  \end{array}
  }
  \end{equation*}
Here, $S\langle q\rangle: \Q\to \Mod{\kk}$ denotes the \emph{stalk functor}, defined by
$$S\langle q\rangle=\Q(q,-)/\mathfrak{r}(q,-)\in\lMod{\Q},$$
and described as follows (see \cite[Lemma 7.10]{HJ22}):
\begin{itemize}
\item for every object $p\in\Q$,
$$S\langle q\rangle(p)=
\begin{cases}
\kk, & \text{if } p=q,\\
0, & \text{otherwise};
\end{cases}$$
\item for every morphism $f\in\Q$,
$$S\langle q\rangle(f) =
\begin{cases}
x\cdot \id_\kk, & \text{if } f=x\cdot\id_q+g\in(\kk\cdot\id_q)\oplus\mathfrak r_q=\Q(q,q),\\
0, & \text{otherwise};
\end{cases}$$
\end{itemize}
The (contravariant) stalk functor $S\{q\}$ is defined dually.
\end{bfhpg}

\begin{bfhpg}[Acyclic objects]\label{acyclic}
Following \cite[Definition 4.1]{HJ22}, let $\lE{\Q,A}$ denote the class of \emph{acyclic objects} (also called \emph{exact objects}) in $\lMod{\Q,A}$ defined as follows:
$$\lE{\Q,A}=\setof{X\in\lMod{\Q,A}}{\pd_{\Q}X^{\natural}<\infty\ \mathrm{or}\ \id_{\Q}X^{\natural}<\infty}.$$
Here, $(-)^{\natural}$ denotes the forgetful functor from $\lMod{\Q,A}$ to $\lMod{\Q}$. It is known that $\lE{\Q,A}$ is a thick subcategory by \cite[Theorem 4.4]{HJ22}, that is, it is closed under direct summands and if two out of three objects in a short exact sequence
$$0 \to X' \to X \to X'' \to 0$$
in $\lMod{\Q,A}$ belong to $\lE{\Q,A}$, then so does the third one.
\end{bfhpg}

The following result, given in \cite[Theorem 7.1]{HJ22}, is used frequently in the paper.

\begin{lem}\label{lem:homology}
Let $X$ be an object in $\lMod{\Q,A}$. Then the following statements are equivalent.
\begin{eqc}
\item $X\in\lE{\Q,A}$.
\item $\Ext[\Q]{i}{S\langle q\rangle}{X}=0$ for each $q\in\Q$ and $i>0$.
\item $\Ext[\Q]{1}{S\langle q\rangle}{X}=0$ for each $q\in\Q$.
\item $\Tor[\Q]{i}{S\{q\}}{X}=0$ for each $q\in\Q$ and $i>0$.
\item $\Tor[\Q]{1}{S\{q\}}{X}=0$ for each $q\in\Q$.
\end{eqc}
\end{lem}

\begin{rmk}\label{rmk:dual-homology}
The dual version of Lemma \ref{lem:homology} also holds for objects $X\in\rMod{\Q,A}$; such an object is acyclic if and only if $\Ext[\Qop]{1}{S\{q\}}{X}=0$ for every $q\in\Q$, or equivalently, $\Tor[\Q]{1}{X}{S\langle q\rangle}=0$ for every $q\in\Q$.
\end{rmk}

\begin{bfhpg}[Complexes]
Let $\sf A$ be an abelian category. In this paper, complexes of objects in $\sf A$, or $\sf A$-complexes, are indexed homologically. That is, the degree $n$ term of an $\sf A$-complex $X$ is denoted $X_n$, and we denote the degree $n$ cycle, cokernel, and homology objects by $Z_n(X)$, $C_n(X)$, and $H_n(X)$, respectively. We denote by $\ch{\sf A}$ the category of $\sf A$-complexes.
\end{bfhpg}

\begin{bfhpg}[Gorenstein projective objects]
Recall from \cite{EnochsGP} that an acyclic complex $T$ of projective objects in $\lMod{\Q,A}$ is called \emph{totally acyclic} if $\Hom_{\Q,A}(T,P)$ is acyclic for each projective object $P$ in $\lMod{\Q,A}$. An object $M$ in $\lMod{\Q,A}$ is called \emph{Gorenstein projective} if it has a complete projective resolution, that is, a totally acyclic complex $T$ of projective objects in $\lMod{\Q,A}$ such that $M\is C_0(T)$; Gorenstein projective $A$-modules may be defined similarly.

It follows from \cite[Corollary 4.8]{DSS17} (or \cite[Theorem B]{DLLM}) that $S\langle q\rangle$ is a Gorenstein projective object in $\lMod{\Q}$, and $S\{q\}$ is a Gorenstein projective object in $\rMod{\Q}$.
\end{bfhpg}

The next definition can be found in \cite[Definitions 2.4 and 2.5]{HJ26}.

\begin{dfn}\label{dfn:good}
A complete projective resolution $P\langle q\rangle$ of $S\langle q\rangle$ is called \emph{good} if it satisfies the following conditions:
\begin{prt}
\item $P\langle q\rangle_0\is\Q(q,-)$ and $P\langle q\rangle_{-1}\is\Q(\mathbb{S}^{-1}q,-)$.
\item The zeroth differential $\partial_0: P\langle q\rangle_0 \to P\langle q\rangle_{-1}$ satisfies $\Hom_{\Q}(S\langle p\rangle, \partial_0)=0$ and $\Hom_{\Q}(\partial_0, S\langle p\rangle)=0$ for each $p\in\Q$.
\end{prt}

Similarly, a complete projective resolution $P\{q\}$ of $S\{q\}$ is called \emph{good} if it satisfies the following conditions:
\begin{prt}
\item[(a')] $P\{q\}_0\is\Q(-,q)$ and $P\{q\}_{-1}\is\Q(-,\mathbb{S}q)$.
\item[(b')] The zeroth differential $\partial_0: P\{q\}_0 \to P\{q\}_{-1}$ satisfies $\Hom_{\Qop}(S\{p\}, \partial_0)=0$ and $\Hom_{\Qop}(\partial_0, S\{p\})=0$ for each $p\in\Q$.
\end{prt}
\end{dfn}

\begin{ipg}\label{finite-generated}
In the paper, we let ${}_\Q\sf prj$ be the full subcategory of finitely generated and
projective objects in $\lMod{\Q}$, that is, objects of the form $\coprod_{q\in V}\Q(q,-)\otimes_\kk M_q$ with $V$ a finite subset of $\Q$ and $M_q$ a finitely generated and projective $\kk$-module; see \cite[Definition 1.1]{HJ26}. By Hom-finiteness and local boundedness of $\Q$, every object in ${}_\Q\mathsf{prj}$ has finite support and finitely generated $\kk$-module values. The subcategory $\sf prj_\Q$ is defined dually.
\end{ipg}

\begin{rmk}\label{rmk:standard}
It follows from \cite[Lemmas 2.7 and 2.9]{HJ26} that, for each $q\in\Q$, the complete projective resolution $P\langle q\rangle$ of $S\langle q\rangle$ can be chosen such that it satisfies Definition \ref{dfn:good}(a), and each of its terms is in ${}_\Q\sf prj$. In this case, $P\{q\}={\sf{\Sigma}^{-1}}\Hom_\kk(P\langle q\rangle,\kk)$ is a complete projective resolution of $S\{q\}$ with each term in $\sf prj_\Q$. Moreover, if $\mathbb{S}q\neq q$, then both $P\langle q\rangle$ and $P\{q\}$ are good.
\end{rmk}

\begin{setup}\label{setup:choosen}
Throughout the paper, for each $q\in\Q$ we choose $P\langle q\rangle$ and $P\{q\}$ to be the complete projective resolutions of $S\langle q\rangle$ and $S\{q\}$ respectively as in Remark \ref{rmk:standard}. In particular, each term of $P\langle q\rangle$ is in ${}_\Q\sf prj$, and each term of $P\{q\}$ is in $\sf prj_\Q$. Thus this choice fits the setup of \cite[Setup 3.1]{HJ26}.
\end{setup}

\begin{bfhpg}[Adjoint triple]\label{adjoint-pair}
For each $q\in\Q$, following \cite{HJ26}, there is an adjoint triple
  \begin{equation*}
  \xymatrix@C=4pc{
    \lMod{\Q,A}
    \ar[r]^-{\nablaq}
    &
    \ch A,
    \ar@/_1.8pc/[l]_-{\DeltaLq}
    \ar@/^1.8pc/[l]^-{\DeltaRq}
  }
  \end{equation*}
where $\nablaq(-)=\Hom_{\Q}(P\langle q\rangle,-)$. By \cite[Lemmas 3.5$(i)$ and 4.5]{HJ26}, all the functors $\DeltaLq$, $\nablaq$ and $\DeltaRq$ are exact, so for each $i\geq0$, there are isomorphisms
\begin{equation*}
\Ext[\ch A]{i}{C}{\nablaq X}\cong\Ext[\Q,A]{i}{\DeltaLq C}{X},
\end{equation*}
and
\begin{equation*}
\Ext[\ch A]{i}{\nablaq X}{C}\cong\Ext[\Q,A]{i}{X}{\DeltaRq C}.
\end{equation*}
\end{bfhpg}

The next result shows that acyclic objects in $\lMod{\Q,A}$ can be detected by these functors in \ref{adjoint-pair}.

\begin{lem}\label{lem:exact}
An object $X$ in $\lMod{\Q,A}$ is acyclic if and only if $\nabla_qX$ is an acyclic $A$-complex for every $q\in\Q$.
\end{lem}
\begin{prf*}
By definition, one gets that $X\in\lE{\Q,A}$ if and only if $X^\natural\in\lE{\Q}$. Moreover, for every $q\in\Q$ there is a natural isomorphism $(\nabla_qX)^\natural\cong\nabla_q(X^\natural)$. Thus the assertion follows from \cite[Lemma 3.12]{HJ26}.
\end{prf*}

\begin{lem}\label{lem:K-cycle}
For every $X\in\lMod{\Q,A}$ and every $q\in\Q$, there is a natural isomorphism $K_q(X)\cong Z_0(\nablaq X)$.
\end{lem}

\begin{prf*}
The complete projective resolution $P\langle q\rangle$ of $S\langle q\rangle$ gives an exact sequence
$$P\langle q\rangle_1\to
P\langle q\rangle_0\to
S\langle q\rangle\to 0.$$
Applying the functor $\Hom_{\Q}(-,X)$, one gets the exact sequence
$$0\to\Hom_{\Q}(S\langle q\rangle,X)
\to\Hom_{\Q}(P\langle q\rangle_0,X)
\to\Hom_{\Q}(P\langle q\rangle_1,X).$$
Thus one has $K_q(X)=\Hom_{\Q}(S\langle q\rangle,X)\cong Z_0(\nablaq X)$.
\end{prf*}

\begin{ipg}\label{right-adjoint}
For each $q\in\Q$, applying the same construction to $\Qop$ gives an adjoint triple
  \begin{equation*}
  \xymatrix@C=4pc{
    \rMod{\Q,A}
    \ar[r]^-{\widetilde{\nabla}_q}
    &
    \ch\Aop,
    \ar@/_1.8pc/[l]_-{{}_q\widetilde{\Delta}}
    \ar@/^1.8pc/[l]^-{\widetilde{\Delta}_q}
  }
  \end{equation*}
where $\widetilde{\nabla}_q(-)=\Hom_{\Qop}(P\{q\},-)$. All the functors ${}_q\widetilde{\Delta}$, $\widetilde{\nabla}_q$ and $\widetilde{\Delta}_q$ are exact.

We also use the adjoint triple $(\DeltaLq^{\kk}, \nablaq^{\kk}, \DeltaRq^{\kk})$ between $\lMod{\Q}$ and $\ch\kk$, and the adjoint triple $({}_q\widetilde{\Delta}^{\kk}, \widetilde{\nabla}_q^{\kk}, \widetilde{\Delta}_q^{\kk})$ between $\rMod{\Q}$ and $\ch\kk$. Here the superscript $\kk$ is used to emphasize that these functors take values in $\kk$-modules.
\end{ipg}

The following two lemmas will be used in the paper.

\begin{lem}\label{lem:serre-hom-comparison}
Let $X$ be a complex of objects in ${}_\Q\sf prj$, $L$ an $A$-module, and $P\in\lMod{\Q,A}$. Then there is a natural isomorphism
$$\Hom_A(\Hom_{\Q}(X,P),L)\cong
\Hom_{\Q^{\mathrm{op}}}(\Hom_{\kk}(X,\kk),
\Hom_A(P,L)\circ\mathbb S^{-1})$$
of complexes of $\kk$-modules. In particular, one has
$$\Hom_A(\nabla_qP,L)\cong{\sf\Sigma}^{-1}\widetilde{\nabla}_q^{\kk}
(\Hom_A(P,L)\circ\mathbb S^{-1}).$$
\end{lem}

\begin{prf*}
It is enough to construct the isomorphism degreewise, naturally in each term. We mention that it suffices to take $X_n=\Q(q,-)\otimes_{\kk}K$ with $K$ a finitely generated and projective $\kk$-module.
Let $K^*=\Hom_{\kk}(K,\kk)$. Then one has
\begin{align*}
\Hom_A(\Hom_{\Q}(X_n,P),L)
&\cong \Hom_A(\Hom_{\kk}(K,P(q)),L)\\
&\cong \Hom_A(K^*\otimes_\kk P(q),L)\\
&\cong \Hom_{\kk}(K^*,\Hom_A(P(q),L)),
\end{align*}
where the second isomorphism holds as $K$ is a finitely generated and projective $\kk$-module. On the other hand, by \cite[Lemma 1.3$(i)$]{HJ26}, one has
$$\Hom_{\kk}(X_n,\kk)\cong \Q(-,\mathbb S q)\otimes_{\kk}K^*.$$
Thus one gets that
\begin{align*}
&\Hom_{\Q^{\mathrm{op}}}(\Hom_{\kk}(X_n,\kk),
\Hom_A(P,L)\circ\mathbb S^{-1})\\
&\qquad\cong
\Hom_{\kk}(K^*,(\Hom_A(P,L)\circ\mathbb S^{-1})(\mathbb S q))\\
&\qquad=
\Hom_{\kk}(K^*,\Hom_A(P(q),L)).
\end{align*}
This gives the first isomorphism of complexes in the statement.

Finally, since $P\langle q\rangle$ is a complex of objects in ${}_\Q\sf prj$ and $P\{q\}={\sf{\Sigma}^{-1}}\Hom_\kk(P\langle q\rangle,\kk)$ by Setup \ref{setup:choosen}, one gets that
$$\Hom_A(\nabla_qP,L)\cong{\sf\Sigma}^{-1}\widetilde{\nabla}_q^{\kk}
(\Hom_A(P,L)\circ\mathbb S^{-1}).$$
Thus the second isomorphism in the statement holds.
\end{prf*}

\begin{lem}\label{lem:mixed-Hom-transfer}
Let $C$ be an $A$-complex and $L$ an $A$-module. Then for each $q\in\Q$, there is a natural isomorphism
$$\widetilde{\Delta}_q^{\kk}\Hom_A(C,L)\cong
\Hom_A(({}_q\Delta{\sf\Sigma} C)\circ\mathbb S^{-1},L).$$
\end{lem}

\begin{prf*}
For every $Y\in\rMod{\Q}$, one has
\begin{align*}
\Hom_{\Q^{\mathrm{op}}}(Y,\widetilde{\Delta}_q^{\kk}\Hom_A(C,L))
&\cong \Hom_{\ch\kk}(\widetilde{\nabla}_q^{\kk}Y,\Hom_A(C,L)) \\
&\cong \Hom_{\ch A}(C,\Hom_{\kk}(\widetilde{\nabla}_q^{\kk}Y,L))\\
&\cong \Hom_{\ch A}({\sf\Sigma} C, \nabla_q(\Hom_{\kk}(Y,L)\circ\mathbb S)) \\
&\cong \Hom_{\Q,A}({}_q\Delta{\sf\Sigma} C, \Hom_{\kk}(Y,L)\circ\mathbb S)\\
&\cong \Hom_{\Q,A}(({}_q\Delta{\sf\Sigma} C)\circ\mathbb S^{-1},\Hom_{\kk}(Y,L))\\
&\cong \Hom_{\Q^{\mathrm{op}}}(Y, \Hom_A(({}_q\Delta{\sf\Sigma} C)\circ\mathbb S^{-1},L)).
\end{align*}
The first and fourth isomorphisms hold as $(\widetilde{\nabla}_q^{\kk}, \widetilde{\Delta}_q^{\kk})$ and $({}_q\Delta, \nabla_q)$ are adjoint pairs; see \ref{right-adjoint} and \ref{adjoint-pair}. The third isomorphism follows from \cite[Lemma 3.6$(iv)$]{HJ26}. The fifth isomorphism is induced by the autoequivalence $\mathbb S$, and the last one is from \cite[Remark 1.4]{HJ26}. All the above isomorphisms are natural in $Y$, so Yoneda's lemma gives the asserted isomorphism.
\end{prf*}

\section{Gorenstein projective criteria for total acyclicity}\label{sec:adjoint-triple}
\noindent
In this section, we give some characterizations for totally acyclic $Q$-shaped diagrams of projective $A$-modules via Gorenstein projective modules. We begin with the next definition given in \cite{HJ26}.

\begin{dfn}\label{dfn:totally-acyclic}
Recall that an object $P\in\lMod{\Q,A}$ is called a \emph{totally acyclic} $\Q$-shaped diagram of projective $A$-modules if $P$ is acyclic with each $P(q)$ projective for $q\in\Q$, and $\Hom_A(P,L)$ is acyclic in $\rMod{\Q}$ for every projective $A$-module $L$.
\end{dfn}

\begin{lem}\label{lem:tacyclic}
If $P$ is a totally acyclic $\Q$-shaped diagram of projective $A$-modules, then $\nabla_qP$ is a totally acyclic complex of projective $A$-modules for each $q\in\Q$.
\end{lem}
\begin{prf*}
Since $P$ is an acyclic object in $\lMod{\Q,A}$ with $P(q)$ projective for each $q\in\Q$, one gets that $\nabla_qP$ is an acyclic complex of $A$-modules by Lemma \ref{lem:exact}, and its terms are projective $A$-modules by \cite[Lemma 3.5$(ii)$]{HJ26}. Let $L$ be a projective $A$-module. Then $\Hom_A(P,L)$ is an acyclic object in $\rMod{\Q}$, and so is $\Hom_A(P,L)\circ\mathbb S^{-1}$. Thus by \cite[Lemma 3.12]{HJ26}, the complex $\widetilde{\nabla}_q^{\kk}(\Hom_A(P,L)\circ\mathbb S^{-1})$ is acyclic. This implies that the complex $\Hom_A(\nabla_qP,L)$ is acyclic by Lemma \ref{lem:serre-hom-comparison}, and so $\nabla_qP$ is a totally acyclic complex of projective
$A$-modules.
\end{prf*}

The next result may be known; we include a proof for the sake of completeness.

\begin{lem}\label{lem:finite-dimension-test}
Let $T$ be an acyclic complex of projective objects in $\rMod{\Q}$ and let $X$ be an object in $\lMod{\Q}$ of finite projective dimension. Then the complex $T\otimes_{\Q}X$ is acyclic.
\end{lem}

\begin{prf*}
We argue by induction on the projective dimension $\pd_{\Q}X$. If $X$ is projective, then it is a direct summand of a coproduct of representables, and so the complex $T\otimes_{\Q}X$ is acyclic as it is a direct summand of a coproduct of evaluation complexes of $T$. Assume that $\pd_{\Q}X=d>0$, and let
$$0\to X'\to P\to X\to 0$$
be an exact sequence in $\lMod{\Q}$ with $P$ projective and $\pd_{\Q}X'=d-1$. Since every term of $T$ is projective, one gets an exact sequence
$$0\to T\otimes_{\Q}X'\to T\otimes_{\Q}P\to T\otimes_{\Q}X\to 0.$$
The first two complexes are acyclic by induction and the projective case, so the third one is acyclic.
\end{prf*}

Now we are in a position to give the main result in this section.

\begin{thm}\label{thm:totally-acyclic-characterization}
Let $P$ be an acyclic object in $\lMod{\Q,A}$ with $P(q)$ projective for each $q\in\Q$. Then the following statements are equivalent.
\begin{eqc}
\item $P$ is a totally acyclic $\Q$-shaped diagram of projective $A$-modules.
\item $K_q(P)$ is Gorenstein projective for every $q\in\Q$.
\item $C_q(P)$ is Gorenstein projective for every $q\in\Q$.
\end{eqc}
\end{thm}

\begin{prf*}
\proofofimp{i}{ii} By Lemma \ref{lem:tacyclic}, one gets that $\nabla_qP$ is a totally acyclic complex of projective $A$-modules for each $q\in\Q$, and so it follows from Lemma \ref{lem:K-cycle} that $K_q(P)\cong Z_0(\nabla_qP)$ is a Gorenstein projective $A$-module.

\proofofimp{ii}{i} Fix $q\in\Q$. Then for $n\in\mathbb Z$, one has
\begin{equation}\tag{\ref{thm:totally-acyclic-characterization}.1}\label{2.3.0}
Z_{-n}(\nablaq P)=Z_{-n}(\Hom_{\Q}(P\langle q\rangle,P))=\Hom_{\Q}(C_n(P\langle q\rangle),P).
\end{equation}
Applying \cite[Construction 7.21(a)]{HJ22} to $C_n(P\langle q\rangle)$, one gets some objects
$$Y^0=C_n(P\langle q\rangle), Y^1, Y^2, \cdots$$
in $\lMod{\Q}$, where $Y^{\ell+1}$ is the third term in the next short exact sequence
\begin{equation}\tag{\ref{thm:totally-acyclic-characterization}.2}\label{2.3.1}
0\to\bigoplus_{p\in \Q}S_pK_p(Y^\ell)\to Y^\ell\to Y^{\ell+1}\to 0
\end{equation}
whenever we have defined $Y^\ell$.

We claim that $Y^N=0$. Indeed, since the pseudo-radical satisfies \(\mathfrak r^N=0\), the definition of the functors $K_p^\ell$ gives $K_p^\ell(Y^0)=Y^0(p)$ for each $\ell\ge N$; see \cite[Definition 7.22]{HJ22}. Thus by \cite[Lemma 7.24(a)]{HJ22}, one has
$$K_p(Y^N)=\operatorname{Coker}(K_p^N(Y^0)\rightarrowtail K_p^{N+1}(Y^0))=0$$
for each $p\in\Q$, and so $Y^N=0$ by \cite[Proposition 7.19]{HJ22}, as desired.

We then collect some facts needed. Since $P\langle q\rangle_n\in{}_\Q\mathsf{prj}$, it has finite support and finitely generated $\kk$-module values (see \ref{finite-generated}), and hence the same holds for $Y^0=C_n(P\langle q\rangle)$. Each $Y^{\ell+1}$ is a quotient of $Y^\ell$, so induction shows that every $Y^\ell$ has finite support and finitely generated $\kk$-module values. On the other hand, since $Y^0=C_n(P\langle q\rangle)$ has projective
$\kk$-module values as $\kk$ is hereditary, it follows from \cite[Proposition 7.20(a)]{HJ22}, using the exact sequence (\ref{2.3.1}), that every $Y^\ell(p)$ and $K_p(Y^\ell)$ are projective over $\kk$. Moreover, there is a containment $K_p(Y^\ell)\subseteq Y^\ell(p)$ by \cite[Proposition 7.18]{HJ22}, so $K_p(Y^\ell)$ is finitely generated as $\kk$ is noetherian, and only finitely many of them are nonzero.

Since $P$ is an acyclic object in $\lMod{\Q,A}$, the object $P^\natural$ has finite projective dimension. By \cite[Corollary 4.8]{DSS17} (see also \cite[Theorem B]{DLLM}), the objects $Y^{\ell+1}$ are Gorenstein projective. Then it follows from
\cite[Theorems 2.2 and 2.4]{HJ22} that $\Ext[\Q]{1}{Y^{\ell+1}}{P^\natural}=0$. Applying $\Hom_{\Q}(-,P)$ to (\ref{2.3.1}), one gets an exact sequence
$$0 \to \Hom_{\Q}(Y^{\ell+1},P)\to\Hom_{\Q}(Y^\ell,P) \to \Hom_{\Q}(\bigoplus_{p\in\Q}S_pK_p(Y^\ell),P) \to 0.$$
Only finitely many $K_p(Y^\ell)$ are nonzero, so the coproduct in the last term is finite. Since $K_p(Y^\ell)$ is finitely generated and projective over $\kk$, one gets that
$$\Hom_{\Q}(S_pK_p(Y^\ell),P)\cong\Hom_{\kk}(K_p(Y^\ell),K_p(P))$$
is a direct summand of a finite coproduct of copies of $K_p(P)$. By assumption, each $K_p(P)$ is Gorenstein projective over $A$. Since the class of Gorenstein projective $A$-modules is closed under finite coproducts, direct summands and extensions (see \cite[Theorem~2.5]{HHl04a}), the last term in the above exact sequence is Gorenstein projective. As $Y^N=0$, backward induction on $\ell$
shows that
$$\Hom_{\Q}(C_n(P\langle q\rangle),P)=\Hom_{\Q}(Y^0,P)$$
is Gorenstein projective, and so every cycle of $\nablaq P$ is Gorenstein projective by (\ref{2.3.0}). On the other hand, since $P$ is an acyclic object in $\lMod{\Q,A}$ with $P(q)$ projective for each $q\in\Q$, the $A$-complex $\nablaq P$ is acyclic by Lemma \ref{lem:exact}, and its terms are projective $A$-modules by \cite[Lemma 3.5$(ii)$]{HJ26}. Thus $\nablaq P$ is a totally acyclic complex of projective $A$-modules.

Let $L$ be a projective $A$-module. Then the complex $\Hom_A(\nabla_qP,L)$ is acyclic, and so by Lemma \ref{lem:serre-hom-comparison}, the complex $\widetilde{\nabla}_q^{\kk}(\Hom_A(P,L)\circ\mathbb S^{-1})$ is acyclic for every $q\in\Q$. Thus by \cite[Lemma 3.12]{HJ26}, the object $\Hom_A(P,L)\circ\mathbb S^{-1}$ in $\rMod{\Q}$ is acyclic, and so $\Hom_A(P,L)$ is acyclic. This yields that $P$ is a totally acyclic
$\Q$-shaped diagram of projective $A$-modules.

\proofofimp{i}{iii} Fix $q\in\Q$. Since $P$ is acyclic in $\lMod{\Q,A}$, one gets that $P^\natural$ has finite projective dimension in $\lMod{\Q}$. Thus by Lemma \ref{lem:finite-dimension-test}, one gets that
$P\{q\}\otimes_{\Q}P^\natural$ is an acyclic complex of $\kk$-modules, which yields that $P\{q\}\otimes_{\Q}P$ is an acyclic complex of $A$-modules. We mention that each $P\{q\}_n\in\sf prj_\Q$ is a direct summand of a finite coproduct of right representables. Then $P\{q\}_n\otimes_{\Q}P$ is a direct summand of a finite coproduct of projective $A$-modules $P(p)$, and so every term of $P\{q\}\otimes_{\Q}P$ is projective over $A$.

Let $L$ be a projective $A$-module. Then one gets that the complex
$$\Hom_A(P\{q\}\otimes_{\Q}P,L)
\cong\Hom_{\Q^{\mathrm{op}}}(P\{q\},\Hom_A(P,L)) =\widetilde{\nabla}_q^{\kk}\Hom_A(P,L)$$
is acyclic by \cite[Lemma 3.12]{HJ26} as $\Hom_A(P,L)$ is an acyclic object in $\rMod{\Q}$. Thus $P\{q\}\otimes_{\Q}P$ is a totally acyclic complex of projective $A$-modules. The nonnegative half of $P\{q\}$ gives a projective resolution
$$P\{q\}_1\to P\{q\}_0\to S\{q\}\to 0$$
of $S\{q\}$. Then one gets that
$$C_q(P)=S\{q\}\otimes_{\Q}P\cong\Coker(P\{q\}_1\otimes_{\Q}P\to
P\{q\}_0\otimes_{\Q}P)=C_0(P\{q\}\otimes_{\Q}P)$$
is Gorenstein projective.

\proofofimp{iii}{i} Let $L$ be a projective $A$-module. For $q\in\Q$, take the projective resolution
$$\cdots\to P\{q\}_1\to P\{q\}_0\to S\{q\}\to 0$$
of $S\{q\}$ given by the nonnegative half of the complete resolution $P\{q\}$, where every $P\{q\}_i$ is in $\sf prj_\Q$; see Setup \ref{setup:choosen}. Since $P$ is in $\lE{\Q,A}$, one has $\Tor[\Q]{i}{S\{q\}}{P}=0$ for each $i>0$ by Lemma \ref{lem:homology}. Thus one gets the exact sequence
\begin{equation}\tag{\ref{thm:totally-acyclic-characterization}.3}\label{2.3.3}
\cdots\to P\{q\}_1\otimes_{\Q}P
\to P\{q\}_0\otimes_{\Q}P
\to C_q(P)\to 0.
\end{equation}
Moreover, each $P\{q\}_i$ is a direct summand of a finite coproduct of right representables, so $P\{q\}_i\otimes_{\Q}P$ is a projective $A$-module as every $P(q)$ is projective by assumption. Thus (\ref{2.3.3}) is a projective resolution of $C_q(P)$ over $A$. This yields the isomorphism
$$\Ext[\Q^{\mathrm{op}}]{i}{S\{q\}}{\Hom_A(P,L)}
\cong\Ext[A]{i}{C_q(P)}{L}$$
for every $i>0$. Since $C_q(P)$ is Gorenstein projective by assumption, the right-hand side vanishes. Since this holds for every $q\in\Q$, one gets that $\Hom_A(P,L)$ is acyclic in $\rMod{\Q}$; see Remark \ref{rmk:dual-homology}. Thus $P$ is a totally acyclic $\Q$-shaped diagram of projective $A$-modules.
\end{prf*}

\section{Stable and defect categories of acyclic Q-shaped diagrams}
\label{sec:tac-stable}
\noindent
In this section, we study the stable categories of acyclic and totally acyclic $Q$-shaped diagrams, introduce the associated defect category, and investigate its relationship with the classical projective defect category of $A$.

\begin{ipg}
Let $\AcQ$ be the subcategory of acyclic objects $P$ in $\lMod{\Q,A}$ with $P(q)$ projective for each $q\in\Q$. We equip this subcategory with the exact structure inherited from $\lMod{\Q,A}$; every short exact sequence whose three terms are in $\AcQ$ is objectwise split as its right-hand term has projective $A$-module values.
\end{ipg}

The next result shows that $\AcQ$ is a Frobenius exact category. We let $\lPrj{\Q,A}$ denote the subcategory of projective objects in $\lMod{\Q,A}$, and for a Frobenius exact category $\sf E$, we write $\mathsf{stab}(\sf E)$ for
its stable category.

\begin{prp}\label{thm:EPrj-Frobenius}
The exact category $\AcQ$ is Frobenius, and its projective-injective objects are precisely the projective objects in
$\lMod{\Q,A}$, that is,
$${\Prj{\AcQ}}={\Inj{\AcQ}}=\lPrj{\Q,A}.$$
Consequently, the stable category $\mathsf{stab}(\AcQ)$ is a triangulated category.
\end{prp}

\begin{prf*}
We first mention that $\lPrj{\Q,A}\subseteq\AcQ$. Indeed, projective objects in $\lMod{\Q,A}$ are acyclic by \cite[Theorems 7.1 and 7.29]{HJ22}, and their evaluations are projective $A$-modules as the evaluation functor $E_q$ has an exact right adjoint $G_q$; see \ref{eva}.

Let $X\in\AcQ$. The counits of the adjunctions $(F_q, E_q)$ induce a canonical epimorphism
$$\varepsilon_X:\mathbf F(X)=\bigoplus_{q\in\Q}F_qX(q)\to X.$$
For each $p\in\Q$, the component $\varepsilon_X(p)$ is split epic, so
\begin{equation}\tag{\ref{thm:EPrj-Frobenius}.1}\label{eq:F-approx}
0\to\mathbf K(X)\to\mathbf F(X)
\xra{\varepsilon_X}X\to0
\end{equation}
is an objectwise split exact sequence; see \cite[Proposition 5.7]{HJ24}. Since $X(q)$ is a projective $A$-module and
$F_q$ preserves projectives by \cite[Lemma 3.11]{HJ22}, the object $\mathbf F(X)$ is projective in $\lMod{\Q,A}$. Thus
$\mathbf K(X)(p)$ is projective for every $p\in\Q$, and $\mathbf K(X)\in\lE{\Q,A}$ as $\lE{\Q,A}$ is a thick subcategory; see \ref{acyclic}. Thus $\mathbf K(X)\in\AcQ$, and so $\AcQ$ has enough projectives.

Dually, the units of the adjunctions $(E_q, G_q)$ give a canonical monomorphism
$$\eta_X: X\to
\mathbf G(X)=\prod_{q\in\Q}G_qX(q).$$
For each $p$, the component $\eta_X(p)$ is split monic, so
\begin{equation}\tag{\ref{thm:EPrj-Frobenius}.2}\label{eq:G-approx}
0\to X\xra{\eta_X}\mathbf G(X)
\to\mathbf C(X)\to0
\end{equation}
is an objectwise split exact sequence; see \cite[Proposition 5.7]{HJ24}. It follows from \cite[Lemma 3.4]{HJ24} that there is an isomorphism $G_q(M)\cong F_{\mathbb S^{-1}q}(M)$ for every $A$-module $M$, so each $G_qX(q)$ is projective.
Moreover, by \cite[Proposition 3.7]{HJ24} one has
$$\bigoplus_{q\in\Q}G_qX(q)\is
\prod_{q\in\Q}G_qX(q),$$
so $\mathbf G(X)$ is projective in $\lMod{\Q,A}$. Then it follows as above that $\mathbf C(X)\in\AcQ$.

It remains to identify the projective-injective objects in $\AcQ$. Let $P\in\lPrj{\Q,A}$. It follows from the proof of \cite[Theorem 7.29(a)]{HJ22} that $P$ may be written as
$$\bigoplus_{q\in\Q}F_q(M_q)\is\prod_{q\in\Q}G_{\mathbb S q}(M_q)$$
with each $M_q$ a projective $A$-module, where the isomorphism holds by \cite[Lemma 3.4 and Proposition 3.7]{HJ24}. For a conflation $\xi$ in $\AcQ$ and every projective $A$-module $M$, one gets that $E_q(\xi)$ is split exact for each $q\in\Q$, so the adjunction
$$\Hom_{\Q,A}(-,G_q(M))\cong\Hom_A(E_q(-),M)$$
yields that the functor $\Hom_{\Q,A}(-,G_q(M))$ sends $\xi$ to an exact sequence. This implies that $G_q(M)$ is injective in $\AcQ$, and so $P$ is injective in $\AcQ$. On the other hand, it is clear that $P$ is projective in $\AcQ$.

Conversely, if $X$ is projective in $\AcQ$, then the deflation in (\ref{eq:F-approx}) splits in $\AcQ$, so $X$ is a direct summand of the projective object $\mathbf F(X)$ and hence belongs to $\lPrj{\Q,A}$. If $X$ is injective in $\AcQ$, then the inflation in (\ref{eq:G-approx}) splits, so $X$ is a direct summand of the projective object $\mathbf G(X)$ and again belongs to $\lPrj{\Q,A}$. Thus one has
$$\Prj{\AcQ}=\Inj{\AcQ}=\lPrj{\Q,A},$$
and so $\AcQ$ is a Frobenius exact category.
\end{prf*}

\begin{ipg}
Let ${\TacQ}$ be the subcategory of totally acyclic $\Q$-shaped diagrams of projective $A$-modules. That is,
$${\TacQ}=\{P\in{\AcQ}\ |\ \Hom_A(P,L) \text{ is acyclic for each }L\in\lPrj A\}.$$
\end{ipg}

In the following, we show that the category $\TacQ$ is also Frobenius.

\begin{lem}\label{lem:projectives-in-Tac}
There is a containment $\lPrj{\Q,A}\subseteq\TacQ$.
\end{lem}

\begin{prf*}
Let $P$ be a projective object in $\lMod{\Q,A}$. Then $P$ is acyclic by \cite[Theorems 7.1 and 7.29]{HJ22} and every $P(q)$ is projective over $A$ as the evaluation functor $E_q$ has an exact right adjoint $G_q$; see \ref{eva}. The functor $S_q$ is exact for each $q\in\Q$, so its left adjoint $C_q$ preserves projective objects. Thus $C_q(P)$ is a projective $A$-module, and hence is a Gorenstein projective $A$-module for each $q\in\Q$. By Theorem \ref{thm:totally-acyclic-characterization}, one has $P\in\TacQ$.
\end{prf*}

\begin{lem}\label{lem:Tac-extension-closed}
The subcategory $\TacQ$ is thick in $\AcQ$.
\end{lem}

\begin{prf*}
The closure under direct summands is immediate from the definition. Consider a conflation
\begin{equation}\tag{\ref{lem:Tac-extension-closed}.1}\label{3.2.1}
0\to P'\to P\to P''\to 0
\end{equation}
in $\AcQ$. It is objectwise split, so for every projective $A$-module $L$, one gets a short exact sequence
\begin{equation}\tag{\ref{lem:Tac-extension-closed}.2}\label{3.2.2}
0\to\Hom_A(P'',L)\to\Hom_A(P,L)\to\Hom_A(P',L)\to 0
\end{equation}
in $\rMod{\Q}$. The subcategory of acyclic objects in $\rMod{\Q}$ is thick; see \ref{acyclic}. Thus, if two of the terms in (\ref{3.2.2}) are exact, then so is the third one. It follows that if two of the three terms in (\ref{3.2.1}) are in $\TacQ$, then so is the third one.
\end{prf*}

\begin{prp}\label{thm:Tac-Frobenius}
The inherited exact structure makes $\TacQ$ a Frobenius exact category, and its projective-injective objects are precisely the projective objects in $\lMod{\Q,A}$, that is,
$$\Prj{\TacQ}=\Inj{\TacQ}=\lPrj{\Q,A}.$$
Consequently, the stable category $\mathsf{stab}(\TacQ)$ is a triangulated category, and the inclusion $\TacQ\hookrightarrow\AcQ$ induces a fully faithful exact triangle functor
$$\mathsf{stab}(\TacQ)\to\mathsf{stab}(\AcQ)$$
whose essential image is thick.
\end{prp}

\begin{prf*}
Let $P$ be in $\TacQ$, and apply (\ref{eq:F-approx}) and (\ref{eq:G-approx}) to $P$. Their middle terms are projective objects in $\lMod{\Q,A}$, hence belong to $\TacQ$ by Lemma \ref{lem:projectives-in-Tac}. Then it follows from Lemma \ref{lem:Tac-extension-closed} that the kernel in (\ref{eq:F-approx}) and the cokernel in (\ref{eq:G-approx}) also belong to $\TacQ$. Consequently, if an object in $\TacQ$ is projective or injective, splitting the
corresponding exact sequence (\ref{eq:F-approx}) or (\ref{eq:G-approx}) yields that
$$\Prj{\TacQ}\subseteq\lPrj{\Q,A}\ \ \text{and}\ \ \Inj{\TacQ}\subseteq\lPrj{\Q,A}.$$
On the other hand, by Proposition \ref{thm:EPrj-Frobenius}, the reverse inclusions hold, so one gets that
$$\Prj{\TacQ}=\Inj{\TacQ}=\lPrj{\Q,A}.$$

Since $\TacQ$ and $\AcQ$ have the same projective-injective objects, the inclusion $\TacQ\hookrightarrow\AcQ$ induces a fully faithful triangle functor
$$\mathsf{stab}(\TacQ)\to\mathsf{stab}(\AcQ).$$
Its essential image is triangulated since $\TacQ$ is closed under extensions and under the suspension and inverse suspension in $\mathsf{stab}(\AcQ)$. To see that the image is thick, suppose that $P\in\AcQ$ is a direct summand of some $T\in\TacQ$ in $\mathsf{stab}(\AcQ)$. Then there exist morphisms 
$$f: P\to T\ \ \text{and}\ \ g: T\to P$$
such that $\mathrm{id}_P-gf$ factors through a
projective-injective object $I\in\lPrj{\Q,A}$. Consequently, $P$ is a direct summand of $T\oplus I$. Since $T\oplus I\in\TacQ$ by Lemma \ref{lem:projectives-in-Tac}, and $\TacQ$ is closed under direct summands by Lemma \ref{lem:Tac-extension-closed}, one gets that $P\in\TacQ$. Thus the essential image is thick.
\end{prf*}

Now let $\mathsf D_{\mathrm{def}}^{\Q}(A)=\mathsf{stab}(\AcQ)\big/\mathsf{stab}(\TacQ)$ be the Verdier quotient, and let $\GP{A}$ denote the subcategory of Gorenstein projective $A$-modules. Then we have the next result.

\begin{thm}\label{cor:defect-zero-criterion}
The following statements are equivalent.
\begin{eqc}
\item $\mathsf D_{\mathrm{def}}^{\Q}(A)=0$.
\item $\AcQ=\TacQ$.
\item For every $P\in\AcQ$ and every $q\in\Q$, one has $C_q(P)\in\GP{A}$.
\item For every $P\in\AcQ$ and every $q\in\Q$, one has $K_q(P)\in\GP{A}$.
\end{eqc}
\end{thm}

\begin{prf*}
The equivalence of \eqclbl{ii}, \eqclbl{iii} and \eqclbl{iv} follows from Theorem \ref{thm:totally-acyclic-characterization}.

The implication \eqclbl{ii} $\Longrightarrow$ \eqclbl{i} is obvious.

\proofofimp{i}{ii} Since $\mathsf{stab}(\TacQ)$ is a thick subcategory of $\mathsf{stab}(\AcQ)$ by Proposition \ref{thm:Tac-Frobenius}, and $\mathsf D_{\mathrm{def}}^{\Q}(A)=0$ by assumption, one gets that every object $P\in\AcQ$ is isomorphic to an object $T\in\TacQ$ in $\mathsf{stab}(\AcQ)$. Then there exist morphisms
$$f: P\to T\ \ \text{and}\ \ g: T\to P$$
such that $\id_P-gf$ factors through a projective-injective object $I\in\lPrj{\Q,A}$. Thus $P$ is a direct summand of $T\oplus I$. Since $T\oplus I$ is in $\TacQ$ by Lemma \ref{lem:projectives-in-Tac}, one gets that $P\in\TacQ$ as $\TacQ$ is closed under direct summands by Lemma \ref{lem:Tac-extension-closed}.
\end{prf*}

\begin{ipg}
Let $K_{\mathrm{ac}}(\Prj A)$ be the homotopy category of acyclic complexes of projective $A$-modules, and let $K_{\mathrm{tac}}(\Prj A)$ be the homotopy category of totally acyclic complexes of projective $A$-modules. Recall that the categories of acyclic and totally acyclic complexes of projective $A$-modules, endowed with the degreewise split exact structures, are Frobenius exact categories. Their projective-injective objects are precisely the projective $A$-complexes, and their stable categories are canonically triangle equivalent to $K_{\mathrm{ac}}(\Prj A)$ and $K_{\mathrm{tac}}(\Prj A)$, respectively.

It is known that the full subcategory $K_{\mathrm{tac}}(\Prj A)$ is thick in $K_{\mathrm{ac}}(\Prj A)$. Define the classical projective defect category of $A$ by
$$\mathsf D_{\mathrm{def}}(A)=K_{\mathrm{ac}}(\Prj A)\big/K_{\mathrm{tac}}(\Prj A).$$
\end{ipg}

Next, we present the relationship between $\mathsf D_{\mathrm{def}}(A)$ and $\mathsf D_{\mathrm{def}}^{\Q}(A)$.

\begin{prp}\label{thm:defect-adjunction}
For every $q\in\Q$, the adjoint pair $({}_q\Delta, \nabla_q)$ induces an adjoint pair of triangle functors
$${}_q\overline{\Delta}:\mathsf D_{\mathrm{def}}(A) \rightleftarrows
\mathsf D_{\mathrm{def}}^{\Q}(A):\overline{\nabla}_q.$$
\end{prp}

\begin{prf*}
We mention that the functors ${}_q\Delta$ and $\nabla_q$ are exact; see \ref{adjoint-pair}. Let $C$ be an acyclic complex of
projective $A$-modules. By \cite[Lemmas 4.8$(i)$ and 4.9]{HJ26}, one has ${}_q\Delta C\in\AcQ$. Conversely, if $P\in\AcQ$, then $\nabla_qP$ is an acyclic complex of projective $A$-modules by Lemma \ref{lem:exact} and \cite[Lemma 3.5$(ii)$]{HJ26}. Thus the adjoint pair $({}_q\Delta, \nabla_q)$ restricts to the acyclic subcategories.

We then turn our attention to the totally acyclic subcategories. Let $C\in K_{\mathrm{tac}}(\Prj A)$ and $L$ be a projective $A$-module. Then the complex $\Hom_A({\sf\Sigma}^{-1}C,L)$ is acyclic. It follows from \cite[Lemma 4.10]{HJ26} that
$\widetilde{\Delta}_q^{\kk}\Hom_A({\sf\Sigma}^{-1}C,L)$ is acyclic in $\rMod{\Q}$, so the object $\Hom_A(({}_q\Delta C)\circ\mathbb S^{-1},L)$ is acyclic in $\rMod{\Q}$ by Lemma \ref{lem:mixed-Hom-transfer}. It follows that $\Hom_A({}_q\Delta C,L)$ is acyclic in $\rMod{\Q}$, and so ${}_q\Delta C\in\TacQ$. Conversely, if $P\in\TacQ$, then by Lemma \ref{lem:tacyclic}, one gets that $\nabla_qP$ is a totally acyclic complex of projective $A$-modules. Thus the functors ${}_q\Delta$ and $\nabla_q$ preserve the corresponding totally acyclic subcategories.

Finally, the functor $\nabla_q$ preserves projectives by \cite[Lemma 3.11(i)]{HJ26}, and so does the functor ${}_q\Delta$, since it is left adjoint to the exact functor $\nabla_q$. Consequently, both functors preserve the projective-injective objects of the corresponding Frobenius exact categories, and hence induce an adjunction on the stable categories $K_{\mathrm{ac}}(\Prj A)$ and $\mathsf{stab}(\AcQ)$. Since these functors also preserve the corresponding totally acyclic subcategories, as proved above, they induce triangle functors on the Verdier quotients. Moreover, the unit and counit of the adjunction on the stable categories descend to the Verdier quotients and still satisfy the triangle identities, so the induced functors form the asserted adjoint pair.
\end{prf*}

In the following, we prove Theorem \ref{thmA} advertised in the introduction.

\begin{prp}\label{cor:defect-detection}
The family of triangle functors $\{\overline{\nabla}_q:\mathsf D_{\mathrm{def}}^{\Q}(A)\to\mathsf D_{\mathrm{def}}(A)\}_{q\in\Q}$ is jointly conservative, equivalently, an object $P$ in $\mathsf D_{\mathrm{def}}^{\Q}(A)$ is zero if and only if $\overline{\nabla}_q(P)=0$ for every $q\in\Q$.
\end{prp}

\begin{prf*}
We only prove the ``if" part. Let $P\in\AcQ$ and suppose that $\overline{\nabla}_q(P)=0$ for every $q\in\Q$, that is, $\nabla_qP$ is in $K_{\mathrm{tac}}(\Prj A)$. Then by Lemma \ref{lem:K-cycle}, one gets that $K_q(P)\cong Z^0(\nabla_qP)$ is Gorenstein projective for each $q\in\Q$. Thus it follows from Theorem \ref{thm:totally-acyclic-characterization} that $P\in\TacQ$, and so one gets that $P=0$ in $\mathsf D_{\mathrm{def}}^{\Q}(A)$.
\end{prf*}

\begin{cor}\label{cor:defect-vanishing-forward}
If $\mathsf D_{\mathrm{def}}(A)=0$, then one has $\mathsf D_{\mathrm{def}}^{\Q}(A)=0$.
\end{cor}

\begin{prf*}
If $\mathsf D_{\mathrm{def}}(A)=0$, then $\overline{\nabla}_q(P)=0$ for each $q\in\Q$ and every object $P$ in $\mathsf D_{\mathrm{def}}^{\Q}(A)$. Thus by Proposition \ref{cor:defect-detection}, one has $\mathsf D_{\mathrm{def}}^{\Q}(A)=0$.
\end{prf*}

\begin{lem}\label{lem:good-detection}
Let $q\in\Q$, and assume that $P\langle q\rangle$ is a good complete projective resolution of $S\langle q\rangle$. Then a $\kk$-complex $D$ is acyclic if and only if $\widetilde{\Delta}_q^{\kk}({\sf\Sigma}^{\ell}D)$ is acyclic in $\rMod{\Q}$ for each $\ell\in\mathbb Z$.
\end{lem}

\begin{prf*}
The complete resolution $P\{q\}={\sf\Sigma}^{-1}\Hom_{\kk}(P\langle q\rangle,\kk)$ of $S\{q\}$ is good by \cite[Lemma 2.9$(ii)$]{HJ26}. Thus the assertion follows from \cite[Lemma 4.10]{HJ26}.
\end{prf*}

\begin{prp}\label{thm:defect-left-conservative}
Let $q\in\Q$, and assume that $P\langle q\rangle$ is a good complete projective resolution of $S\langle q\rangle$. Then the triangle functor ${}_q\overline{\Delta}: \mathsf D_{\mathrm{def}}(A)\to\mathsf D_{\mathrm{def}}^{\Q}(A)$ is conservative, equivalently, an object $C$ in $\mathsf D_{\mathrm{def}}(A)$ is zero if and only if ${}_q\overline{\Delta}(C)=0$.
\end{prp}

\begin{prf*}
Let $C\in K_{\mathrm{ac}}(\Prj A)$ and suppose that ${}_q\overline{\Delta}(C)=0$. Then one gets that ${}_q\overline{\Delta}({\sf\Sigma}^{\ell}C)=0$ for all $\ell\in\mathbb Z$, which yields that the stable class of ${}_q\Delta{\sf\Sigma}^{\ell}C$ is in $\mathsf{stab}(\TacQ)$ for every $\ell\in\mathbb Z$. Equivalently, ${}_q\Delta{\sf\Sigma}^{\ell}C$ is stably isomorphic to an object $T_\ell\in\TacQ$, which implies that ${}_q\Delta{\sf\Sigma}^{\ell}C$ is a direct summand of $T_\ell\oplus I$ for some projective-injective
$I\in\lPrj{\Q,A}$. By Lemma \ref{lem:projectives-in-Tac} and the fact that $\TacQ$ is closed under direct summands, one gets that ${}_q\Delta{\sf\Sigma}^{\ell}C\in\TacQ$. Let $L$ be a projective $A$-module. Then it follows from Lemma \ref{lem:mixed-Hom-transfer} that
$$\widetilde{\Delta}_q^{\kk}\Hom_A({\sf\Sigma}^{\ell-1}C,L)\cong \Hom_A(({}_q\Delta{\sf\Sigma}^{\ell}C)\circ\mathbb S^{-1},L)$$
is acyclic in $\rMod{\Q}$ for every $\ell\in\mathbb Z$. Thus by Lemma \ref{lem:good-detection}, one gets that the complex $\Hom_A(C,L)$ is acyclic, and so $C$ is totally acyclic. Thus one gets that $C=0$ in $\mathsf D_{\mathrm{def}}(A)$. This implies that the functor ${}_q\overline{\Delta}$ is conservative.
\end{prf*}

\begin{cor}\label{cor:defect-vanishing-equivalence}
Assume that there exists $q\in\Q$ for which $P\langle q\rangle$ is good. Then
$\mathsf D_{\mathrm{def}}(A)=0$ if and only if $\mathsf D_{\mathrm{def}}^{\Q}(A)=0$.
\end{cor}

\begin{prf*}
The ``only if" part is Corollary \ref{cor:defect-vanishing-forward}. For the ``if" part, we assume that $\mathsf D_{\mathrm{def}}^{\Q}(A)=0$ and choose $q\in\Q$ with $P\langle q\rangle$ good. Then one has
${}_q\overline{\Delta}(C)=0$ for every $C\in\mathsf D_{\mathrm{def}}(A)$, and so by
Proposition \ref{thm:defect-left-conservative}, one gets that $C=0$ in $\mathsf D_{\mathrm{def}}(A)$. Thus $\mathsf D_{\mathrm{def}}(A)=0$.
\end{prf*}

\begin{thm}\label{thm:classical-Q-defect-equivalence}
Consider the following conditions:
\begin{eqc}
\item $K_{\mathrm{ac}}(\Prj A)=K_{\mathrm{tac}}(\Prj A)$.
\item $\AcQ=\TacQ$.
\end{eqc}
Then \eqclbl{i}$\Longrightarrow$\eqclbl{ii} holds. If, in addition, $\mathbb S q\neq q$ for some $q\in\Q$, then the two conditions are equivalent.
\end{thm}

\begin{prf*}
By Theorem \ref{cor:defect-zero-criterion}, one gets that $\AcQ=\TacQ$ if and only if $\mathsf D_{\mathrm{def}}^{\Q}(A)=0$. Thus the implication \eqclbl{i}$\Longrightarrow$\eqclbl{ii} holds by Corollary \ref{cor:defect-vanishing-forward}. Conversely, choose $q\in\Q$ such that $\mathbb S q\neq q$. Then $P\langle q\rangle$ is good by \cite[Lemma 2.7$(iii)$]{HJ26}. Thus by Corollary \ref{cor:defect-vanishing-equivalence}, the implication \eqclbl{ii}$\Longrightarrow$\eqclbl{i} holds.
\end{prf*}

\begin{exa}\label{exa:converse-fails}
The implication \eqclbl{ii}$\Longrightarrow$\eqclbl{i} in Theorem \ref{thm:classical-Q-defect-equivalence} does not hold in general without an additional hypothesis. Indeed, let $\kk$ be a field and let $\Q$ be the $\kk$-linear category with one object $q$ and $\Q(q,q)=\kk$. Then $\mathbb S=\id_{\Q}$ and the pseudo-radical of $\Q$ is zero. Moreover, every $\Q$-module is exact as it is just a $\kk$-module. Thus for every $\kk$-algebra $A$, one has
$$\AcQ=\Prj A=\TacQ.$$

Now take $A=\kk[x,y]/(x,y)^2$. Then $A$ is a commutative artinian local ring with socle $(x,y)$ of $\kk$-dimension two, and so it is not Gorenstein as an artinian local Gorenstein ring has one-dimensional socle. By \cite[Theorem 19.5.25]{CFH24}, one gets that
$$K_{\mathrm{ac}}(\Prj A)\neq K_{\mathrm{tac}}(\Prj A).$$
Thus the condition \eqclbl{ii} of Theorem \ref{thm:classical-Q-defect-equivalence} does not imply the condition \eqclbl{i} in general. In particular, $\mathbb S q=q$ in this example.
\end{exa}

For a commutative noetherian ring $A$, we use the term \emph{Gorenstein} in the usual local sense; that is, $A_{\mathfrak p}$ is Gorenstein for every $\mathfrak p\in\operatorname{Spec}A$. The equivalence of \eqclbl{i} and \eqclbl{iii} in the next result was given by Holm and J{\o}rgensen; see \cite[Theorem B]{HJ26}.

\begin{cor}\label{cor:Gorenstein-Q-defect}
Assume that $A$ is commutative noetherian, and $\mathbb S q\neq q$ for some $q\in\Q$. Then the following conditions are equivalent.
\begin{eqc}
\item $A$ is Gorenstein.
\item $K_{\mathrm{ac}}(\Prj A)=K_{\mathrm{tac}}(\Prj A)$.
\item $\AcQ=\TacQ$.
\end{eqc}
\end{cor}

\begin{prf*}
Under the assumption that $A$ is commutative noetherian, the equivalence \eqclbl{i}$\Longleftrightarrow$\eqclbl{ii} holds by e.g. \cite[Theorem 19.5.25]{CFH24}, and the equivalence \eqclbl{ii}$\Longleftrightarrow$\eqclbl{iii} follows from Theorem \ref{thm:classical-Q-defect-equivalence}.
\end{prf*}

\begin{cor}\label{cor:finite-id-Q-defect}
Assume that $\id_A L<\infty$ for every projective $A$-module $L$. Then one has $\AcQ=\TacQ$.
\end{cor}

\begin{prf*}
Let $C$ be an acyclic complex of projective $A$-modules and $L$ a projective $A$-module, and let $\id_A L=d<\infty$. For every $n$, consider the exact sequence
$$0\to Z_nC\to C_n\to Z_{n-1}C \to 0.$$
By dimension shifting, one has
$$\Ext[A]{1}{Z_{n-1}C}{L}
\cong\Ext[A]{d+1}{Z_{n-d-1}C}{L}=0.$$
This implies that the complex $\Hom_A(C,L)$ is acyclic. Thus one gets that
$$K_{\mathrm{ac}}(\Prj A)=K_{\mathrm{tac}}(\Prj A),$$
and so by Theorem \ref{thm:classical-Q-defect-equivalence}, one has $\AcQ=\TacQ$.
\end{prf*}

\section{Applications to N-complexes and N-stable categories}
\noindent
Let $\Q_{N\text{-}\mathrm{cpx}}$ be the $\kk$-linear category given by the quiver
$$\cdots\to 2\to 1\to 0\to -1\to -2\to\cdots$$
with the relations that every composition of $N$ consecutive arrows is zero for $N\geq2$. Then the category $\Q_{N\text{-}\mathrm{cpx}}$ satisfies the hypotheses in Setup \ref{setup:Q}; see \cite{HJ22}. The $\Q_{N\text{-}\mathrm{cpx}}$-shaped diagrams are precisely $N$-complexes, and so one has $\lMod{\Q_{N\text{-}\mathrm{cpx}},A}={\sf Ch}^N(A)$, the category of $N$-complexes of $A$-modules. Moreover, acyclic objects in $\lMod{\Q_{N\text{-}\mathrm{cpx}},A}$ are precisely the $N$-acyclic $N$-complexes of $A$-modules; see \cite[Section 3.5]{HJ24b}. In this section, we apply the results of the previous two sections to
$N$-complexes and obtain a strengthening of the $N$-stable equivalence of Brightbill and Miemietz \cite[Theorem 4.12]{BM24}.
\begin{ipg}
An $N$-complex $P$ is written as $P=(P_n,d_n:P_n\to P_{n-1})_{n\in\mathbb Z}$ with
$$d_{n-N+1}\cdots d_{n-1}d_n=0.$$
For $1\leq r<N$, we write
$$d_n^{\,r}=d_{n-r+1}\cdots d_{n-1}d_n: P_n\to P_{n-r}.$$
Following the convention of \cite{BM24}, we let $Z_n^r(P)=\Ker d_n^{\,r}$, $B_n^r(P)=\operatorname{Im}d_{n+N-r}^{\,N-r}$, and $H_n^r(P)=Z_n^r(P)/B_n^r(P)$. Recall that $P$ is called $N$-\emph{acyclic} if $H_n^r(P)=0$ for all $n\in\mathbb Z$ and $1\leq r<N$.

{\bf Throughout the section}, let $\mathsf{Acy}^N(\Prj A)$ denote the category of $N$-acyclic $N$-complexes of projective $A$-modules, and let $\mathsf{Tac}^N(\Prj A)$ be the full subcategory of $\mathsf{Acy}^N(\Prj A)$ consisting of those $P$ for which $\Hom_A(P,L)$ is $N$-acyclic for every projective $A$-module $L$; the $N$-complexes in $\mathsf{Tac}^N(\Prj A)$ are called \emph{totally $N$-acyclic}. Let $K^N(\Prj A)$ be the homotopy category of $N$-complexes of projective $A$-modules, and write $K^N_{\mathrm{ac}}(\Prj A)$ and $K^N_{\mathrm{tac}}(\Prj A)$ for the full subcategories of $K^N(\Prj A)$ consisting of $N$-acyclic and totally $N$-acyclic $N$-complexes.
\end{ipg}

\begin{rmk}\label{rmk:N-total-terminology}
Our terminology is adapted to the $\Q$-shaped notion of total acyclicity used throughout the paper. It should not be confused with the terminology of Bahiraei, Hafezi and Nematbakhsh in \cite[Section 4]{BHN16}, where an $N$-totally acyclic complex of finitely generated projectives is defined by two-sided homotopy orthogonality to bounded projective $N$-complexes. We will not need to identify these two notions here.
\end{rmk}

We first make explicit the functors $K_n$ and $C_n$ in \ref{stalk} in the present special case.

\begin{lem}\label{lem:N-stalks}
Let $P$ be an $N$-complex and let $n\in\mathbb Z$. Then one has $K_n(P)\cong Z_n^1(P)$
and $C_n(P)\cong\Coker(d_{n+1}: P_{n+1}\to P_n)$. Moreover, if $P$ is $N$-acyclic, then one has
$$C_n(P)\cong Z_{n-N+1}^1(P).$$
\end{lem}

\begin{prf*}
The left stalk $S\langle n\rangle$ is concentrated at the vertex $n$. Hence a morphism $S\langle n\rangle\to P$ is determined by an element of $P_n$ annihilated by $d_n$, which yields that $$K_n(P)\cong\Ker d_n=Z_n^1(P).$$ 

Dually, tensoring with the right stalk $S\{n\}$ gives $C_n(P)\cong P_n/\operatorname{Im}d_{n+1}$.
If $P$ is $N$-acyclic, then one has $\operatorname{Im}d_{n+1}=B_n^{N-1}(P)=Z_n^{N-1}(P)$. Thus
$$C_n(P)\cong P_n/Z_n^{N-1}(P) \cong\operatorname{Im}d_n^{\,N-1}
=B_{n-N+1}^{1}(P) =Z_{n-N+1}^{1}(P),$$
as required.
\end{prf*}

Now Theorem \ref{thm:totally-acyclic-characterization} gives the following result, which will be used later to extend the $N$-stable equivalence studied by Brightbill and Miemietz \cite{BM24} beyond the Gorenstein ambient case.

\begin{thm}\label{cor:N-cycles-GP}
Let $P$ be an $N$-acyclic $N$-complex of projective $A$-modules. Then the following conditions are equivalent.
\begin{eqc}
\item $P$ is totally $N$-acyclic.
\item $K_n(P)$ is Gorenstein projective for every $n\in\mathbb Z$.
\item $C_n(P)$ is Gorenstein projective for every $n\in\mathbb Z$.
\item $Z_n^1(P)$ is Gorenstein projective for every $n\in\mathbb Z$.
\item $Z_n^r(P)$ is Gorenstein projective for every $n\in\mathbb Z$ and all
$1\leq r<N$.
\end{eqc}
\end{thm}

\begin{prf*}
Since exact objects in $\lMod{\Q_{N\text{-}\mathrm{cpx}},A}$ are precisely $N$-acyclic $N$-complexes (see \cite[3.5]{HJ24b}), the equivalence of \eqclbl{i}, \eqclbl{ii} and \eqclbl{iii} is Theorem \ref{thm:totally-acyclic-characterization} applied to $\Q_{N\text{-}\mathrm{cpx}}$. By Lemma \ref{lem:N-stalks}, \eqclbl{ii} is equivalent to \eqclbl{iv}.

The implication \eqclbl{v} $\Rightarrow$ \eqclbl{iv} is clear. It remains to prove \eqclbl{iv} $\Rightarrow$ \eqclbl{v}. For $1\leq r<N-1$, the $N$-acyclicity of $P$ yields an exact sequence
\begin{equation}\tag{\ref{cor:N-cycles-GP}.1}\label{4.4.1}
0\to Z_n^r(P)\to Z_n^{r+1}(P)\xra{d_n^{\,r}}
Z_{n-r}^1(P)\to 0.
\end{equation}
Indeed, the kernel of $d_n^{\,r}$ is clearly $Z_n^r(P)$. To prove surjectivity, let $y\in Z_{n-r}^1(P)$. Since $P$ is $N$-acyclic, one has $$Z_{n-r}^1(P)=B_{n-r}^1(P)=\operatorname{Im}d_{n-r+N-1}^{\,N-1}.$$
Thus there exists $z\in P_{n-r+N-1}$ such that $y=d_{n-r+N-1}^{\,N-1}(z)=d_n^{\,r}\,d_{n-r+N-1}^{\,N-r-1}(z)$. Let $x=d_{n-r+N-1}^{\,N-r-1}(z)\in P_n$. Then one has $y=d_n^{\,r}(x)$. Moreover,
$$0=d_{n-r}(y)=d_n^{\,r+1}(x),$$
so $x\in Z_n^{r+1}(P)$. Thus $d_n^{\,r}$ is surjective. Since the class of Gorenstein projective $A$-modules is closed under extensions by \cite[Theorem 2.5]{HHl04a}, condition \eqclbl{iv} implies \eqclbl{v} by induction on $r$.
\end{prf*}

Propositions \ref{thm:EPrj-Frobenius} and \ref{thm:Tac-Frobenius} in the previous section also specialize directly.

\begin{prp}\label{cor:N-Frobenius}
With the inherited exact structures, both $\mathsf{Acy}^N(\Prj A)$ and $\mathsf{Tac}^N(\Prj A)$ are Frobenius exact categories; the projective-injective objects are precisely the contractible $N$-complexes of projective $A$-modules. Moreover, there are canonical triangle equivalences
$$\mathsf{stab}(\mathsf{Acy}^N(\Prj A))
\simeq K^N_{\mathrm{ac}}(\Prj A)
\ \ \text{and}\ \
\mathsf{stab}(\mathsf{Tac}^N(\Prj A))
\simeq K^N_{\mathrm{tac}}(\Prj A),$$
and $K^N_{\mathrm{tac}}(\Prj A)$ is a thick subcategory of $K^N_{\mathrm{ac}}(\Prj A)$.
\end{prp}

\begin{prf*}
Applying Propositions \ref{thm:EPrj-Frobenius} and \ref{thm:Tac-Frobenius} to $\Q_{N\text{-}\mathrm{cpx}}$, one gets that $\mathsf{Acy}^N(\Prj A)$ and $\mathsf{Tac}^N(\Prj A)$ are Frobenius exact categories whose projective-injective objects are precisely the projective objects in ${\sf Ch}^N(A)$. These are exactly the contractible $N$-complexes of projective $A$-modules. Indeed, every projective object in ${\sf Ch}^N(A)$ is also projective for the degreewise split exact structure, and hence is contractible by \cite[Theorem 2.1]{IKM17}. Conversely, every contractible $N$-complex of projective $A$-modules is a direct sum of elementary contractible $N$-complexes with projective components, each of which is projective in ${\sf Ch}^N(A)$; see \cite[Theorem 2.1]{IKM17}.

By \cite[Theorem 2.6]{IKM17}, null-homotopic morphisms of $N$-complexes are precisely those factoring through projective-injective objects. For morphisms between $N$-complexes of projective modules, the standard factorization of a null-homotopic morphism may be chosen through a contractible $N$-complex of projective modules. Thus passing to the stable categories is exactly the same as passing to $N$-homotopy classes. Consequently, one has
$$\mathsf{stab}(\mathsf{Acy}^N(\Prj A))
\simeq K^N_{\mathrm{ac}}(\Prj A)
\ \ \text{and}\ \
\mathsf{stab}(\mathsf{Tac}^N(\Prj A))
\simeq K^N_{\mathrm{tac}}(\Prj A).$$

Finally, under the above equivalences, one gets that $K^N_{\mathrm{tac}}(\Prj A)$ is a thick subcategory of $K^N_{\mathrm{ac}}(\Prj A)$ by Proposition \ref{thm:Tac-Frobenius}.
\end{prf*}

\begin{ipg}\label{monomorphism}
Recall that the exact category $\GP{A}$ of Gorenstein projective $A$-modules, with the exact structure inherited from $\lMod{A}$, is Frobenius; its projective-injective objects are precisely the projective $A$-modules. Let $\operatorname{MMor}^{N-2}(\GP{A})$ denote the exact category of chains of $N-2$ admissible monomorphisms in $\GP{A}$. By \cite[Theorem 3.12]{BM24}, this is again Frobenius. Following \cite[Definition 3.13]{BM24}, define
$$\mathsf{stab}^N(\GP{A}):=\mathsf{stab}(\operatorname{MMor}^{N-2}(\GP{A})).$$
For $N=2$, this is the stable category $\mathsf{stab}(\GP{A})$.

For later use, we mention that if $P$ is totally $N$-acyclic, then all its generalized cycles belong to $\GP A$ by Theorem \ref{cor:N-cycles-GP}. Thus for every $n\in\mathbb Z$ and $1\leq r<N$ there is an exact sequence
\begin{equation}\tag{\ref{monomorphism}.1}\label{4.7.1}
0\to Z_n^r(P)\to P_n\xra{d_n^{\,r}}Z_{n-r}^{N-r}(P)\to 0.
\end{equation}
Together with (\ref{4.4.1}), these sequences show that all monomorphisms and epimorphisms between generalized cycles used below are admissible in $\GP A$.
\end{ipg}

We next obtain a totally acyclic extension of the $N$-stable equivalence of Brightbill and Miemietz \cite[Theorem 4.12]{BM24}, without assuming that the coefficient category is Gorenstein. Notice that our convention $Z_n^r$ corresponds to the notation $Z_r^n$ of \cite{BM24}.

\begin{thm}\label{thm:N-stable-total}
For every $N\geq2$, the functor
$$F_N:\mathsf{Tac}^N(\Prj A)\to\operatorname{MMor}^{N-2}(\GP{A})$$
defined by $F_N(P)=(Z_0^1(P)\hookrightarrow Z_0^2(P)\hookrightarrow\cdots
\hookrightarrow Z_0^{N-1}(P))$ induces a triangle equivalence
$$\overline F_N: K^N_{\mathrm{tac}}(\Prj A)\xra{\sim}
\mathsf{stab}^N(\GP{A}).$$
\end{thm}

\begin{prf*}
Let $P\in\mathsf{Tac}^N(\Prj A)$. Then by Theorem \ref{cor:N-cycles-GP}, one has $Z_n^r(P)\in\GP{A}$ for every $n\in\mathbb Z$ and $1\leq r<N$. Thus the canonical inclusions
$$Z_0^1(P)\hookrightarrow Z_0^2(P)\hookrightarrow\cdots
\hookrightarrow Z_0^{N-1}(P)$$
are admissible monomorphisms in $\GP{A}$, so $F_N$ is well defined. Let
$$0\to P\to Q\to T\to 0$$
be a (degreewise split) conflation in $\mathsf{Tac}^N(\Prj A)$. For every $1\leq r<N$, the argument in the proof of
\cite[Proposition 4.7]{BM24} yields an exact sequence
$$0\to Z_0^r(P)\to Z_0^r(Q)\to Z_0^r(T)\to 0$$
in $\GP{A}$. Thus $F_N$ is exact. We next prove that $F_N$ is full and essentially surjective.

Let $P$ and $Q$ be in $\mathsf{Tac}^N(\Prj A)$, and let $\alpha: F_N(P)\to F_N(Q)$ be a morphism in $\operatorname{MMor}^{N-2}(\GP{A})$. The last component
$$Z_0^{N-1}(P)\to Z_0^{N-1}(Q)\hookrightarrow Q_0$$
extends across the admissible monomorphism $Z_0^{N-1}(P)\hookrightarrow P_0$ as $Q_0$ is projective-injective in $\GP{A}$. This gives a map $f_0:P_0\to Q_0$ whose restriction to $Z_0^{N-1}(P)$ is the composite
$$Z_0^{N-1}(P)\xra{\alpha_{N-1}} Z_0^{N-1}(Q)\hookrightarrow Q_0.$$
Since the components of $\alpha$ commute with the canonical inclusions $Z_0^r\hookrightarrow Z_0^{r+1}$, it follows that, for every $1\leq r<N$, the restriction of $f_0$ to $Z_0^r(P)$ is the composite $$Z_0^r(P)\xra{\alpha_r}Z_0^r(Q)\hookrightarrow Q_0.$$

By Theorem \ref{cor:N-cycles-GP}, all generalized cycles of $P$ and $Q$ are in $\GP{A}$, while all terms $P_n$ and $Q_n$ are projective-injective in $\GP{A}$. Consequently, the recursive argument in \cite[Proposition 4.3]{BM24} applies verbatim in the Frobenius exact category $\GP{A}$. Starting from $f_0$, it produces maps $f_n:P_n\to Q_n$ for all $n\in\mathbb Z$ satisfying
$$d_n^Qf_n=f_{n-1}d_n^P$$
and preserving all generalized cycles. Thus they form a morphism $f: P\to Q$ with $F_N(f)=\alpha$. This implies that $F_N$ is full.

Let $X=(X_1\rightarrowtail\cdots\rightarrowtail X_{N-1}) \in\operatorname{MMor}^{N-2}(\GP{A})$. The construction in
\cite[Proposition 4.5]{BM24}, which uses only the Frobenius exact structure, extends $X$ by successive pushouts and pullbacks to a bi-infinite $N$-acyclic array in $\GP{A}$. It therefore yields an $N$-acyclic $N$-complex $P$ of projective $A$-modules, such that $Z_n^r(P)\in\GP{A}$ for all $n\in\mathbb Z$ and each $1\leq r<N$, and $F_N(P)\cong X$. By Theorem \ref{cor:N-cycles-GP}, one gets that $P$ is totally $N$-acyclic, so $F_N$ is essentially surjective.

By Proposition \ref{cor:N-Frobenius}, the stable category of $\mathsf{Tac}^N(\Prj A)$ is canonically $K^N_{\mathrm{tac}}(\Prj A)$. If $I$ is a contractible $N$-complex of projective $A$-modules, then each $Z_0^r(I)$ is projective and each inclusion $Z_0^r(I)\hookrightarrow Z_0^{r+1}(I)$ is split. Thus $F_N(I)$ is projective-injective in $\operatorname{MMor}^{N-2}(\GP{A})$ by \cite[Proposition 3.9]{BM24}. So $F_N$ induces a triangle functor
$$\overline F_N: K^N_{\mathrm{tac}}(\Prj A)\to\mathsf{stab}^N(\GP{A}).$$

It remains to prove faithfulness. Let $f: P\to Q$ be a morphism in $\mathsf{Tac}^N(\Prj A)$ such that $\overline F_N([f])=0$. Then $F_N(f)$ factors through a projective-injective object in
$\operatorname{MMor}^{N-2}(\GP{A})$. The construction in \cite[Lemma 4.10]{BM24} produces maps
$$h_j: P_j\to Q_{j+N-1}$$
for $2-N\leq j\leq0$, forming a seed for $f$ at degree $0$. The recursive argument of \cite[Lemma 4.11]{BM24}, which again uses only the relevant conflations in $\GP{A}$ and the projective-injectivity of the terms, extends this seed in both directions to maps
$$h_n: P_n\to Q_{n+N-1}$$
for $n\in\mathbb Z$ satisfying the $N$-homotopy identities in every degree. Thus $f$ is null-homotopic, and so $\overline F_N$ is faithful.

Since fullness and essential surjectivity pass to the stable quotients, $\overline F_N$ is full and essentially surjective as well. Thus $\overline F_N$ is a triangle equivalence.
\end{prf*}

\begin{rmk}
We mention that for $N=2$, Theorem \ref{thm:N-stable-total} recovers the classical equivalence $K_{\mathrm{tac}}(\Prj A)\simeq\mathsf{stab}(\GP{A}).$
\end{rmk}

The next corollary recovers the module category case of
\cite[Theorem 4.12]{BM24}.

\begin{cor}\label{cor:BM-recovery}
Assume that $A$ is Iwanaga-Gorenstein. Then there exists a triangle equivalence
$K^N_{\mathrm{ac}}(\Prj A)=K^N_{\mathrm{tac}}(\Prj A)\xra{\sim}\mathsf{stab}^N(\GP{A}).$
\end{cor}

\begin{prf*}
Since $A$ is Iwanaga-Gorenstein, one gets that $\lMod{A}$ is a Gorenstein abelian category. It follows from \cite[Proposition 4.1]{BM24} that every generalized cycle of an $N$-acyclic $N$-complex of projective $A$-modules is Gorenstein projective. Thus by Theorem \ref{cor:N-cycles-GP}, one gets that every such $N$-complex is totally $N$-acyclic, that is, $K^N_{\mathrm{ac}}(\Prj A)=K^N_{\mathrm{tac}}(\Prj A)$. Then the assertion holds by Theorem \ref{thm:N-stable-total}.
\end{prf*}

By Proposition \ref{cor:N-Frobenius}, the $\Q$-shaped defect category from the previous section specializes to
$$\mathsf D_{\mathrm{def}}^{N}(A)
:=\mathsf D_{\mathrm{def}}^{\Q_{N\text{-}\mathrm{cpx}}}(A)
\simeq K^N_{\mathrm{ac}}(\Prj A)\big/K^N_{\mathrm{tac}}(\Prj A).$$

\begin{cor}\label{cor:N-defect-Nstable}
There is a fully faithful triangle functor
$$\iota_N: \mathsf{stab}^N(\GP{A}) \to K^N_{\mathrm{ac}}(\Prj A)$$
whose essential image is $K^N_{\mathrm{tac}}(\Prj A)$. Consequently, there is a triangle equivalence
$$\mathsf D_{\mathrm{def}}^{N}(A) \simeq K^N_{\mathrm{ac}}(\Prj A)\big/\operatorname{EssIm}(\iota_N).$$
\end{cor}

\begin{prf*}
Fix a triangle quasi-inverse $G_N: \mathsf{stab}^N(\GP A) \to K^N_{\mathrm{tac}}(\Prj A)$ of the equivalence $\overline F_N$ in Theorem \ref{thm:N-stable-total}, and let $j: K^N_{\mathrm{tac}}(\Prj A) \hookrightarrow K^N_{\mathrm{ac}}(\Prj A)$
be the canonical fully faithful inclusion. Let $\iota_N=jG_N$. Then $\iota_N$ is a fully faithful triangle functor, and since $G_N$ is essentially surjective, its essential image is precisely $K^N_{\mathrm{tac}}(\Prj A)$. Consequently, one has $\mathsf D_{\mathrm{def}}^{N}(A)\simeq K^N_{\mathrm{ac}}(\Prj A)\big/\operatorname{EssIm}(\iota_N)$.
\end{prf*}

Finally, Proposition \ref{thm:defect-adjunction} specializes as follows. For a commutative noetherian ring, we use the term Gorenstein in the usual local sense.

\begin{cor}\label{cor:N-complex-defect}
The adjoint pair $(\DeltaLq, \nablaq)$ between ${\sf Ch}^N(A)$ and $\ch A$ induces, for every $q\in\Q_{N\text{-}\mathrm{cpx}}$, an adjoint pair of triangle functors
$${}_q\overline{\Delta}:\mathsf D_{\mathrm{def}}(A)
\rightleftarrows \mathsf D_{\mathrm{def}}^{N}(A):\overline{\nabla}_q,$$
with ${}_q\overline{\Delta}$ conservative. Consequently, $\mathsf D_{\mathrm{def}}(A)=0$ if and only if $\mathsf D_{\mathrm{def}}^{N}(A)=0$. If, in addition, $A$ is commutative noetherian, then these conditions are equivalent to $A$ being Gorenstein.
\end{cor}

\begin{prf*}
For the category $\Q_{N\text{-}\mathrm{cpx}}$, one has $\mathbb S q\neq q$ for every object $q$, so $P\langle q\rangle$ is good by \cite[Lemma 2.7$(iii)$]{HJ26}. Now the adjunction follows from Proposition \ref{thm:defect-adjunction}, while the conservativity of ${}_q\overline{\Delta}$ follows from Lemma \ref{thm:defect-left-conservative}. The vanishing equivalence is Corollary \ref{cor:defect-vanishing-equivalence}. Finally, when $A$ is commutative noetherian, $\mathsf D_{\mathrm{def}}(A)=0$ is equivalent to $A$ being Gorenstein; see \cite[Theorem 19.5.25]{CFH24}.
\end{prf*}

\section*{Acknowledgment}
\noindent
We thank Zhenxing Di, Henrik Holm, Peter J{\o}rgensen, Liping Li and Yajun Ma for valuable comments on a draft of this paper.


\end{document}